\documentclass[12pt]{article}
\usepackage{amsthm, amsmath,amssymb,amsfonts}
\usepackage{graphicx}
\usepackage{color}

\usepackage{psfig}
\usepackage{graphics}
\usepackage{amsmath}
\usepackage{amssymb}
\usepackage{epsfig}
\usepackage{epstopdf}
\usepackage{color}

\newtheorem{theorem}{Theorem}[section]
\newtheorem{proposition}[theorem]{Proposition}
\newtheorem{corollary}[theorem]{Corollary}
\newtheorem{lemma}[theorem]{Lemma}
\newtheorem{remark}[theorem]{Remark}
\newtheorem{definition}[theorem]{Definition}
\newtheorem{example}[theorem]{Example}

\def\para{\vspace{2mm}}
\def\cP{{\mathcal P}}

\def\gcd{{\rm gcd}}
\def\lcm{{\rm lcm}}

\def\Resultant{{\rm resultant}}

\def\ord{{\rm ord}}

\def\deg{{\rm deg}}

\def\resultant{{\rm resultant}}

\begin{document}
\title{Asymptotes of Space Curves\thanks{The author S. P\'erez-D\'{\i}az is member of the Research Group ASYNACS (Ref. CCEE2011/R34)}}
\author{Angel Blasco and Sonia P\'erez-D\'{\i}az\\
Departamento de F\'{\i}sica y Matem\'aticas \\
        Universidad de Alcal\'a \\
      E-28871 Madrid, Spain  \\
angel.blasco@uah.es, sonia.perez@uah.es
}
\date{}          
\maketitle

\begin{abstract}
In this paper, we generalize the results presented in \cite{BlascoPerezII} for the case of  real algebraic  space curves. More precisely,   given an algebraic space curve $\cal C$  implicitly defined, we show how to compute   the {\it generalized
asymptotes}. In addition, we show how to deal with this problem for the case of a given curve $\cal C$ parametrically defined. The approaches are  based on the notion of {\it approaching  curves} introduced in \cite{paper1}.
\end{abstract}



 \section{Introduction}

In the first part of the paper (Sections 2, 3 and 4), we consider $\cal C$  an irreducible real algebraic space curve over the field of complex numbers $\Bbb C$  implicitly defined
 by two irreducible polynomials $f_1(x_1, x_2, x_3), f_2(x_1, x_2, x_3)\in {\Bbb R}[x_1, x_2, x_3]$. That is,   we  work
over $\Bbb C$, but $\cal C$ has infinitely many points in the affine plane over the field of real numbers $\Bbb R$. Since every irreducible real curve has a real defining polynomial, we assume that $\cal C$ is defined by  irreducible polynomials in ${\Bbb R}[x_1, x_2, x_3]$  (see Chapter 7 in \cite{SWP}).
\para

In  the second part of the paper (Section 5), we are given an irreducible real algebraic space curve  $\cal C$  defined by a parametrization of the form  $\cP(s)=(p_1(s),p_2(s),p_3(s))$, where $p_i(s)=p_{i1}(s)/p(s),\,i=1,2,3$. Similarly as above, since every real curve  can be  parametrized over $\Bbb R$, we assume that $\cP(s)\in {\Bbb R}(s)^3$  (see Chapter 7 in \cite{SWP}).

\para

In both cases, the assumption of reality is included because of the nature of the
problem, but the theory can be similarly developed for the case of
complex non-real curves.

 \para

Under these conditions, we deal with the problem of computing  the asymptotes of the infinity branches of $\cal C$.   Intuitively speaking, the asymptotes of some branch of an
algebraic curve reflect the status of this branch at the points with
sufficiently large coordinates. In analytic geometry, an asymptote
of a curve is a line such that the distance between the curve and
the line approaches zero as they tend to infinity.  In some
contexts, such as algebraic geometry, an asymptote is defined as a
line which is tangent to a curve at infinity. Thus, the problem of  computing  the asymptotes   is  very important in the study of real  algebraic
curves since  asymptotes contain much of the information about the behavior of the curves in the large.





\para

Determining the linear asymptotes of an  algebraic curve  is a topic considered in many text-books on
analysis (see e.g \cite{maxwell}). In \cite{jovan}, it is presented a  simple method for obtaining the linear asymptotes of a curve defined by an irreducible polynomial, with emphasis on
second order polynomials.  In \cite{Zeng},  an algorithm for computing all the linear asymptotes of a real plane
algebraic curve   implicitly defined, is obtained. In \cite{RSS}, it is  briefly studied the linear asymptotes  of space curves. In particular, it is proved  how the tangents at the simple points at infinity of the curve (i.e. non-singular points at infinity) are related with the asymptotes. 

 \para

 However, an algebraic   curve may have more general curves than lines that describe the status of a branch at the points
with sufficiently large coordinates. This motivates the interest  in analyzing and computing these {\it  generalized asymptotes}. Intuitively speaking, a curve $\widetilde{{\cal C}}$ is a {\it  generalized asymptote} (or \textit{g-asymptote}) of another curve $\cal C$ at some infinity branch $B\subset {\cal C}$ if the distance between $\widetilde{{\cal C}}$ and $B$  tends to zero as they tend to infinity, and $\cal C$  can not be approached by a new curve of lower degree at $B$ (see \cite{BlascoPerezII}).

\para

 A deeply elaborated theory in this sense is developed by the authors in \cite{BlascoPerezII}. In that paper, a method for computing all the g-asymptotes of a real plane algebraic curve $\cal C$
implicitly  defined by an irreducible polynomial $f(x_1,x_2)\in {\Bbb R}[x_1,x_2]$ is presented. The approach is based on the notion of approaching curves introduced  in \cite{paper1}.

\para

In this paper,   we generalize these  results, and we  present an algorithm for computing
g-asymptotes of a real  algebraic space curve $\cal C$ implicitly defined
 by two irreducible polynomials $f_1(x_1,x_2,x_3), f_2(x_1,x_2,x_3)\in {\Bbb R}[x_1,x_2,x_3]$.  In addition, we also show how to compute the g-asymptotes if the given curve is defined
 parametrically. This parametric approach can be easily generalized for parametric plane curves and in general, for a rational parametrization of a curve in the $n$-dimensional space.

 \para

The structure of the paper is as follows. In Section 2, we present the notation and we  generalize some previous results developed in \cite{paper1}. In particular, we
characterize   whether two  implicit  algebraic space
curves approach each other at the infinity, and we present a method
to compare the asymptotic behavior  of two space curves (i.e., the behavior
  at the infinity).  In Section 3, we show the relation between infinity branches of plane curves and infinity branches of space
curves. More precisely, we  obtain the infinity branches of a given space curve $\cal C$ from the infinity branches
of a certain plane curve obtained by projecting $\cal C$
along some ``valid projection direction''.

\para

The study of approaching curves and convergent branches leads to the
notions of \textit{perfect curve} (a curve of degree $d$ that cannot
be approached by any curve of degree less than $d$) and
\textit{g-asymptote} (a perfect curve that approaches another curve
at an infinity branch). These concepts are introduced in Section 4.
In this section, we obtain an algorithm that computes a
g-asymptote for each infinity branch of a given curve.  Section 5 is devoted to the computation of g-asymptotes for a given parametric space curve. We remark that the method presented
in this section is easily applicable to parametric plane curves and in general, for  rational parametrizations of curves in the $n$-dimensional space.

 \section{Notation and terminology}\label{S-notation}

In this section, we present some notions and terminology that will
be used throughout the paper. In particular, we need some previous
results concerning local parametrizations and Puiseux series.  For
further details see \cite{Alon92},  \cite{paper1}, \cite{Duval89},
Section 2.5 in \cite{SWP}, \cite{Stad00},  and Chapter 4 (Section 2)
in \cite{Walker}.

\para

We denote by ${\Bbb C}[[t]]$ the domain of {\em formal power series}
in the indeterminate $t$ with coefficients in the field ${\Bbb C}$,
i.e. the set of all sums of the form $\sum_{i=0}^\infty a_it^i$,
 $a_i \in {\Bbb C}$. The quotient field of ${\Bbb C}[[t]]$ is
called the field of {\em formal Laurent series}, and it is denoted by
${\Bbb C}((t))$. It is well known that every non-zero formal Laurent
series $A \in {\Bbb C}((t))$ can be written in the form
$A(t) = t^k\cdot(a_0+a_1t+a_2t^2+\cdots), {\rm\ where\ }
    a_0\not= 0 {\rm\ and\ } k\in \Bbb Z.$
  In addition, the field
${{\Bbb C}\ll t\gg}:= \bigcup_{n=1}^\infty {\Bbb C}((t^{1/n}))$ is
called the field of {\em  formal Puiseux series}. Note that Puiseux
series are power series of the form
$$\varphi(t)=m+a_1t^{N_1/N}+a_2t^{N_2/N}+a_3t^{N_3/N} +
\cdots\in{\Bbb C}\ll t\gg,\quad a_i\not=0,\, \forall i\in {\Bbb
N},$$
 where  $N, N_i\in {\Bbb N},\,\,i\geq 1$,  and $0<N_1<N_2<\cdots$. The natural number $N$ is known as
 {\em  the ramification index} of the series. We denote it as $\nu(\varphi)$ (see \cite{Duval89}).

\para

The {\em order} of a non-zero (Puiseux or Laurent) series $\varphi$ is the
smallest exponent of a term with non-vanishing coefficient in $\varphi$.
We denote  it  by $\ord(\varphi)$. We let the order of 0 be $\infty$.

\para

The most important property of Puiseux series is given by Puiseux's Theorem, which states that if $\Bbb K$ is an algebraically closed field, then the field $K\ll x\gg$ is algebraically closed  (see Theorems 2.77 and 2.78 in \cite{SWP}). A proof of Puiseux's Theorem can be given constructively by the
Newton Polygon Method (see e.g. Section 2.5 in \cite{SWP}).

\para

Let ${\cal C}\in {\Bbb C}^3$ be an irreducible space curve defined
by two polynomials $f_1(x_1,x_2,x_3),f_2(x_1,x_2,x_3)\in {\Bbb R}[x_1,x_2,x_3]$. We assume that ${\cal C}$ is not planar (for planar space curves, one may apply the results in \cite{BlascoPerezII}).

\para

We note that   we  work
over   $\Bbb C$, but we assume that the curve has infinitely many points in the affine plane over $\Bbb R$ and then,  $\cal C$ has a real defining polynomial (see Chapter 7 in \cite{SWP}). We recall that the assumption of reality is included because of the nature of the
problem, but the theory developed in this paper can be applied  for the case of
complex non-real curves.

\para

 Let
${\cal C}^*$  be the corresponding projective curve  defined by the
homogeneous polynomials $F_1(x_1,x_2,x_3,x_4),$  $F_2(x_1,x_2,x_3,x_4)\in {\Bbb
R}[x_1,x_2,x_3,x_4]$. Furthermore, let $P=(1:m_2:m_3:0),\,m_2,m_3\in {\Bbb
C},$ be an infinity point of ${\cal C}^*$.

\para

In addition,  we consider the curve
defined implicitly by the polynomials $g_i(x_2,x_3,x_4):=F_i(1,x_2,x_3,x_4)\in
{\Bbb R}[x_2,x_3,x_4]$, for $i=1,2$. Observe that $g_i(p)=0,$ where
$p=(m_2,m_3,0)$.  Let $I\in {\Bbb R}(x_4)[x_2,x_3]$ be the ideal generated by $g_i(x_2,x_3,x_4),\,\,i=1,2$ in the ring ${\Bbb R}(x_4)[x_2,x_3]$. Since ${\cal C}$ is not contained in some
hyperplane $x_4 =c,\,c\in {\Bbb C}$, we have that  $x_4$ is not algebraic over $\Bbb R$. Under this assumption,
the ideal $I$ (i.e. the system of equations $g_1 =  g_2 = 0$) has only finitely many
solutions in the $3$-dimensional affine space over the algebraic  closure of ${\Bbb R}(x_4)$ (which is contained in ${{\Bbb C}\ll x_4\gg}$). Then, there are finitely many pairs of Puiseux
series $(\varphi_{2}(t),\varphi_{3}(t)) \in {{\Bbb C}\ll t\gg}^2$
such that $g_i(\varphi_{2}(t),\varphi_{3}(t),t) = 0,\,i=1,2$. Each
of the pairs $(\varphi_{2}(t),\varphi_{3}(t))$ is a solution of the
system, and $\varphi_{2}(t)$ and $\varphi_{3}(t)$ converge in a
neighborhood of $t= 0$.

\para

It is important to remark that if
$\varphi(t):=(\varphi_{2}(t),\varphi_{3}(t))$ is a solution of the
system, then
$\sigma_{\epsilon}(\varphi)(t):=(\sigma_{\epsilon}(\varphi_{2})(t),\sigma_{\epsilon}(\varphi_{3})(t))$
is another  solution of the system, where
$$\sigma_{\epsilon}(\varphi_{k})(t)=\sum_{i\geq 0}a_{i,k}\epsilon^{\lambda_{i,k}} t^{N_{i,k}/N_k},\,\,N_k,\,N_{i,k}\in {\Bbb N},\,\,0<N_{1,k}<N_{2,k}<\cdots,$$
$N:=\lcm(N_2,N_3)$, $\lambda_{i,k}:=N_{i,k}N/N_k\in {\Bbb N}$, and
$\epsilon^N=1$ (see \cite{Alon92}). We refer to these solutions as
the {\em conjugates} of $\varphi$. The set of all (distinct)
conjugates of $\varphi$ is called the  {\em conjugacy class} of
$\varphi$, and the number of different conjugates of $r$ is $N$.  We
denote the natural number $N$  as  $\nu(\varphi)$.

\para

Under these conditions and reasoning as in \cite{paper1}, we get
that there exists   $M \in {\Bbb R}^+$ such that  for $i\in \{1,2\}$,
$$F_i(1:\varphi_{2}(t):\varphi_{3}(t):t)=g_i(\varphi_{2}(t),\varphi_{3}(t),t)=0,\,\, \mbox{for\, $t\in {\Bbb C}$\, and
$|t|<M$},$$
where
$$\varphi_{k}(t)= \sum_{i\geq 0}a_{i,k}t^{N_{i,k}/N_k},\,\,N_k,\,N_{i,k}\in {\Bbb N},\,\,0<N_{1,k}<N_{2,k}<\cdots.$$
This implies that
$$F_i(t^{-1}:t^{-1}\varphi_{2}(t):t^{-1}\varphi_{3}(t):1)=f_i(t^{-1},t^{-1}\varphi_{2}(t),t^{-1}\varphi_{3}(t))=0,$$ for
$t\in {\Bbb C}$  and  $0<|t|<M$.

\para

\noindent
Now, we set $t^{-1}=z$, and   we obtain
that  for $i\in \{1,2\}$,
$$f_i(z,r_{2}(z),r_{3}(z))=0,\quad \mbox{$z\in {\Bbb C}$\, and
$|z|>M^{-1}$,\qquad where}$$
$$r_{k}(z)=z\varphi_{k}(z^{-1})=m_kz+a_{1,k}z^{1-N_{1,k}/N_k}+a_{2,k}z^{1-N_{2,k}/N_k}+a_{3,k}z^{1-N_{3,k}/N_k} + \cdots,$$
$a_{j,k}\not=0$, $N_k,N_{j,k}\in {\Bbb N},\,\,j=1,\ldots$, and $0<N_{1,k}<N_{2,k}<\cdots$.

\para

\para

  Since $\nu(\varphi)=N$, we get that there are $N$ different series in its conjugacy
class. Let  $\varphi_{j,k},\,j=1,\ldots,N$ be these series,  and $r_{j,k}(z)=z\varphi_{j,k}(z^{-1})=$
\begin{equation}\label{Eq-conjugates}
m_kz+a_{1,k}c_j^{\lambda_{1,k}}
z^{1-N_{1,k}/N_k}+a_{2,k}c_j^{\lambda_{2,k}}
z^{1-N_{2,k}/N_k}+a_{3,k}c_j^{\lambda_{3,k}} z^{1-N_{3,k}/N_k} +
\cdots\end{equation} where $N:=\lcm(N_2,N_3)$,
$\lambda_{i,k}:=N_{i,k}N/N_k\in {\Bbb N}$, and $c_1,\ldots,c_N$ are
the $N$ complex roots of $x^N=1$. Now we are ready to introduce the
notion of infinity branch. The following definitions and results generalize those presented in \cite{paper1} for algebraic plane curves.

\para

\begin{definition}\label{D-infinitybranch}
 An {\em
infinity branch of a space curve ${\cal C}$}  associated to the
infinity point $P=(1:m_2:m_3:0),\,m_2,m_3\in {\Bbb C}$, is  a set
$ B=\bigcup_{j=1}^N L_j$, where
$L_j=\{(z,r_{j,2}(z),r_{j,3}(z))\in {\Bbb C}^3: \,z\in {\Bbb
C},\,|z|>M\}$,\,  $M\in {\Bbb R}^+$, and the series $r_{j,2}$ and
$r_{j,3}$ are given by (\ref{Eq-conjugates}).
 The subsets $L_1,\ldots,L_N$ are
called the {\em leaves} of the infinity branch $B$.
 \end{definition}

\begin{remark} \label{R-conjugation}
An infinity branch  is uniquely determined  from one leaf, up to conjugation. That is, if
$B=\bigcup_{j=1}^N L_i$, where
$L_i=\{(z,r_{i,2}(z),r_{i,3}(z))\in {\Bbb C}^3: \,z\in {\Bbb
C},\,|z|>M\}$,  and
$$r_{i,k}(z)=z\varphi_{i,k}(z^{-1})=
m_kz+a_{1,k}
z^{1-N_{1,k}/N_k}+a_{2,k}
z^{1-N_{2,k}/N_k}+a_{3,k} z^{1-N_{3,k}/N_k} +
\cdots$$
then $r_{j,k}=r_{i,k},\,j=1,\ldots,N$, up to conjugation; i.e. $r_{j,k}(z)=z\varphi_{j,k}(z^{-1})=$
$$
m_kz+a_{1,k}c_j^{\lambda_{1,k}}
z^{1-N_{1,k}/N_k}+a_{2,k}c_j^{\lambda_{2,k}}
z^{1-N_{2,k}/N_k}+a_{3,k}c_j^{\lambda_{3,k}} z^{1-N_{3,k}/N_k} +
\cdots$$ $N, N_{i,k}\in\mathbb{N}$, $\lambda_{i,k}:=N_{i,k}N/N_k\in
{\Bbb N},\,k=2,3,$ and $c_j^N=1,\,\,j=1,\ldots,N$.
\end{remark}

\para

\begin{remark} \label{R-infinitypoint} Observe that the above approach and Definition \ref{D-infinitybranch} is presented for infinity points of the
form $(1: m_2: m_3: 0)$. For the infinity points $(0: m_2: m_3: 0)$, with $m_2\not=0$ or $m_3\not=0$, we reason similarly but we dehomogenize w.r.t $x_2$ (if $m_2\not=0$) or $x_3$ (if $m_3\not=0$). More precisely, we distinguish two different cases:
\begin{enumerate}
\item If $(0: m_2: m_3: 0),\,m_2\not=0$ is an infinity point of the given space curve ${\cal C}$, we consider  the curve
defined  by the polynomials $g_i(x_1,x_3,x_4):=F_i(x_1,1,x_3,x_4)\in
{\Bbb R}[x_1,x_3,x_4],\,i=1,2$, and we reason as above. We get that   an
infinity branch of   ${\cal C}$   associated to the
infinity point $P=(0:m_2:m_3:0),\,m_2\not=0$, is  a set
$ B=\bigcup_{j=1}^N L_j$, where
$L_j=\{(r_{j,1}(z),z,r_{j,3}(z))\in {\Bbb C}^3: \,z\in {\Bbb
C},\,|z|>M\}$,\,  $M\in {\Bbb R}^+$. 
\item If $(0: m_2: m_3: 0),\,m_3\not=0$ is an infinity point of the given space curve ${\cal C}$, we consider  the curve
defined  by the polynomials $g_i(x_1,x_2,x_4):=F_i(x_1,x_2,1,x_4)\in
{\Bbb R}[x_1,x_2,x_4],\,i=1,2$, and we reason as above. We get that   an
infinity branch of   ${\cal C}$   associated to the
infinity point $P=(0:m_2:m_3:0),\,m_3\not=0$, is  a set
$B=\bigcup_{j=1}^N L_j$, where
$L_j=\{(r_{j,1}(z),r_{j,2}(z),z)\in {\Bbb C}^3: \,z\in {\Bbb
C},\,|z|>M\}$,\,  $M\in {\Bbb R}^+$.
\end{enumerate}
Additionally, instead of working with this type of branches, if the  space curve $\cal C$
 has infinity points of the form $(0 : m_2 : m_3 : 0)$, one may consider a linear change of coordinates. Thus, in the following, we may  assume w.l.o.g that the given algebraic space curve $\cal C$ only  has infinity points of the form $(1 : m_2 : m_3 : 0)$. More details on this type of branches are given in \cite{paper1}.
\end{remark}

\para

In the following, we introduce the notions of convergent branches
and approaching curves. Intuitively speaking, two infinity branches
converge if they get closer  as they tend to infinity. This concept
will allow us to analyze whether two space curves approach each other and it generalizes the notion introduced for the plane case (see \cite{paper1}).

\begin{definition}\label{D-distance0}
Two
infinity branches, $B$ and $\overline{B}$, are convergent if there
exist two leaves   $L=\{(z,r_2(z), r_3(z))\in {\Bbb C}^3:\,z\in {\Bbb C},\,|z|>M\} \subset B$ and
$\overline{L}=\{(z,\overline{r}_2(z),\overline{r}_3(z))\in {\Bbb C}^3:\,z\in {\Bbb
C},\,|z|>\overline{M}\}\subset \overline{B}$  such that  $$\lim_{z\rightarrow\infty} d(({r}_2(z),{r}_3(z)), (\overline{r}_2(z),\overline{r}_3(z)))=0.$$ In this case, we say that the leaves $L$ and $\overline{L}$ converge.
\end{definition}

\begin{remark} \label{R-distance0}
\begin{enumerate}
\item In Definition \ref{D-distance0}, we consider any distance $d(u,v),\,u,v\in {\Bbb C}^2$. Taking into account that all norms are equivalent in ${\Bbb C}^2$, we easily get that
$\lim_{z\rightarrow\infty} d(({r}_2(z),{r}_3(z)), (\overline{r}_2(z),\overline{r}_3(z)))=0$ if and only if $\lim_{z\rightarrow\infty} (r_i(z)-\overline{r}_i(z))=0,\,i=2,3.$
\item Two convergent infinity branches are associated to
the same infinity point (see Remark 4.5 in \cite{paper1}).
\end{enumerate}
\end{remark}

\para

In the following lemma, we characterize the convergence of two given infinity branches. This result is obtained similarly as in the case of plane curves and thus, we omit the proof (see Lemma 4.2, and Proposition 4.6 in \cite{paper1}).

\begin{lemma}\label{L-DistVertical} The following statements hold:
\begin{itemize}\item Two  leaves $L=\{(z,r_2(z), r_3(z))\in {\Bbb C}^3:\,z\in {\Bbb C},\,|z|>M\}$ and
$\overline{L}=\{(z,\overline{r}_2(z),\overline{r}_3(z))\in {\Bbb C}^3:\,z\in {\Bbb
C},\,|z|>\overline{M}\}$ are convergent if and only if the terms
with non negative exponent in the series $r_i(z)$ and
$\overline{r}_i(z)$ are the same, for $i=2,3$.
\item  Two infinity branches $B$ and $\overline{B}$ are convergent if and
only if for each leaf $L\subset B$ there exists a leaf $\overline{L}\subset
\overline{B}$ convergent with $L$, and reciprocally.\end{itemize}\end{lemma}

\para

In Definition \ref{D-distance1}, we introduce the notion of {\it approaching curves} that is,   {\it curves that approach each other}. For this purpose, we recall that given an algebraic space curve ${\cal C}$ over $\Bbb C$ and a point $p\in {\Bbb
C}^3$,   {\em the distance from $p$ to ${\cal C}$} is defined as
$d(p,{\cal C})=\min\{d(p,q):q\in{\cal C}\}.$

\para

\begin{definition}\label{D-distance1}
Let ${\cal C} $ be an algebraic space curve over ${\Bbb C}$ with an
infinity branch $B$. We say that a  curve ${\overline{{\cal C}}}$
{\em approaches} ${\cal C}$ at its infinity branch $B$ if there
exists one leaf $L=\{(z,r_2(z), r_3(z))\in {\Bbb C}^3:\,z\in {\Bbb
C},\,|z|>M\}\subset B$ such that
$\lim_{z\rightarrow\infty}d((z,r_2(z), r_3(z)),\overline{\cal C})=0.$
\end{definition}

\para

In the following, we state some important results concerning two curves that approach each other. These results can be proved similarly as in the case of plane curves (see Lemma 3.6, Theorem 4.11, Remark 4.12 and  Corollary 4.13 in \cite{paper1}).

\begin{theorem}\label{T-curvas-aprox}
Let ${\cal C}$
be a space algebraic curve over $\Bbb C$ with an infinity branch
$B$. A space algebraic curve ${\overline{{\cal C}}}$ approaches
${\cal C}$ at $B$ if and only if ${\overline{{\cal C}}}$ has an
infinity branch, $\overline{B}$, such that $B$ and $\overline{B}$
are convergent.
\end{theorem}

 \para

\begin{remark}\label{R-approaching-curves}
\begin{enumerate}
\item Note that ${\overline{{\cal C}}}$ approaches ${\cal C}$ at
some infinity branch $B$ if and only if ${\cal C}$ approaches
${\overline{{\cal C}}}$ at some infinity branch $\overline{B}$. In the following, we
say that ${\cal C}$ and ${\overline{{\cal C}}}$ approach each other
or that they are {\em approaching curves}.
\item Two approaching curves have
  a common infinity point.
\item  ${\overline{{\cal C}}}$ {\em
approaches} ${\cal C}$ at an infinity branch $B$ if and only if for every leaf $L=\{(z,r_2(z), r_3(z))\in {\Bbb C}^3:\,z\in {\Bbb
C},\,|z|>M\}\subset B$, it holds that
$\lim_{z\rightarrow\infty}d((z,r_2(z), r_3(z)),\overline{\cal
C})=0$.
\end{enumerate}
\end{remark}

\para

\begin{corollary}\label{C-approaching-curves}
Let $\cal C$ be an  algebraic space curve with an infinity branch
$B$. Let ${\overline{{\cal C}}}_1$ and ${\overline{{\cal C}}}_2$ be two
different curves that approach $\cal C$ at $B$. Then:
\begin{enumerate}
\item  ${\overline{{\cal C}}}_i$ has an infinity branch
$\overline{B_i}$ that converges with $B$, for $i=1,2$.
\item $\overline{B_1}$ and $\overline{B_2}$ are convergent. Then,   ${\overline{{\cal C}}}_1$
and ${\overline{{\cal C}}}_2$ approach each other.
\end{enumerate}\end{corollary}

\para

For the sake of simplicity, and taking into account that an infinity branch $B$  is uniquely
determined  from one leaf, up to conjugation (see statement 1 in
Remark \ref{R-conjugation}),  we identify an infinity branch by
just one of its leaves. Hence, in the following
 $$B=\{(z,r_2(z), r_3(z))\in {\Bbb C}^3:\,z\in {\Bbb
C},\,|z|>M\},\qquad M\in {\Bbb R}^+$$
will stand for the infinity branch whose leaves are obtained by
conjugation on
$$r_{k}(z)=
m_kz+a_{1,k}
z^{1-N_{1,k}/N_k}+a_{2,k}
z^{1-N_{2,k}/N_k}+a_{3,k} z^{1-N_{3,k}/N_k} +
\cdots,$$
$a_{i,k}\not=0,\, \forall i\in {\Bbb N},\,i\geq 1,$ $N_k, N_{i,k}\in {\Bbb N},\,\,k=2,3$, and $0<N_{1,k}<N_{2,k}<\cdots$ for $k=2,3$.
Observe that   the results stated above hold
for any leaf of $B$. In addition, we also will show that the results obtained in the following sections hold for any leaf (see statement $3$ in Remark \ref{R-algoritmo}).

\section{Computation of infinity branches}

Let $\cal C$ be an irreducible algebraic space curve defined by the polynomials
$f_1(x_1,x_2,x_3),\,f_2(x_1,x_2,x_3)\in {\Bbb R}[x_1, x_2, x_3]$. In  \cite{Bajaj}, it is proved that there exists a plane curve, say ${\cal C}^p$, which is
birationally related with ${\cal C}$. That is, there exists a
birational correspondence between the points of ${\cal C}^p$ and the
points of ${\cal C}$. Furthermore, it is  shown that ${\cal C}^p$ can
always be obtained by projecting ${\cal C}$ along some ``valid
projection direction''.

\para

In the following we assume that the $x_3$-axis is a valid projection
direction (otherwise, we apply a linear change of
coordinates). Let ${\cal C}^p$ be the projection of $\cal C$ along
the $x_3$-axis, and let $f^p(x_1,x_2)\in {\Bbb R}[x_1, x_2]$ be the implicit polynomial defining ${\cal C}^p$. In  \cite{Bajaj}, it is shown  how to construct a birational mapping $h(x_1,x_2)=h_1(x_1,x_2)/h_2(x_1,x_2)$ such that $(x_1,x_2,x_3)\in{\cal C}$ if and only if $(x_1,x_2)\in{\cal C}^p$ and $x_3=h(x_1,x_2)$. We refer to $h(x_1, x_2)$ as the {\em lift function}, since we can obtain the
points of the space curve ${\cal C}$ by applying $h$ to the points of the
plane projected curve ${\cal C}^p$. In addition, note that $x_3=h(x_1,x_2)$ if and only if $h_1(x_1,x_2)-h_2(x_1,x_2)x_3=0$. Thus, ${\cal C}$ can be
implicitly defined by the polynomials $f^p(x_1,x_2)$ and
$f_3(x_1,x_2,x_3)=h_1(x_1,x_2)-h_2(x_1,x_2)x_3.$

\para

In Theorem \ref{T-space-projected-branches}, we study the relation between the infinity branches
of $\cal C$ and ${\cal C}^p$. The idea is to use the lift function
$h$ to obtain the infinity branches of the space curve $\cal C$ from
the infinity branches of  the plane curve ${\cal C}^p$. An efficient method to
compute the infinity branches of a plane curve is presented in
\cite{paper1}.

\begin{theorem}\label{T-space-projected-branches}
$B^p=\{(z,r_2(z))\in {\Bbb C}^2:\,z\in {\Bbb
C},\,|z|>M^p\}$  is an infinity branch of ${\cal C}^p$
for some $M^p\in {\Bbb R}^+$ iff there exists a series $r_3(z)=z\varphi_3(1/z)$,
$\varphi_3(z)\in{{\Bbb C}\ll z\gg}$, such that
$B=\{(z,r_2(z), r_3(z))\in {\Bbb C}^3:\,z\in {\Bbb
C},\,|z|>M\}$ is an infinity branch of $\cal C$
for some $M\in{\Bbb R}^+$.
\end{theorem}

\noindent\textbf{Proof:} Clearly, if $B$ is an infinity branch of $\cal C$, then $B^p$ is an infinity branch of ${\cal C}^p$. Conversely, let
$B^p=\{(z,r_2(z))\in {\Bbb C}^2:\,z\in {\Bbb
C},\,|z|>M^p\}$ be an infinity branch  of  ${\cal C}^p$, and we look for a
series $r_3(z)=z\varphi_3(1/z)$,
$\varphi_3(z)\in{{\Bbb C}\ll z\gg}$, such
that $B=\{(z,r_2(z), r_3(z))\in {\Bbb C}^3:\,z\in {\Bbb
C},\,|z|>M\}$ is an infinity branch of $\cal
C$. Note that, from the discussion above, we can get it as
$r_3(z)=h(z,r_2(z))$. However, we need to prove that
$r_3(z)=z\varphi_3(1/z)$ for some Puiseux series $\varphi_3(z)$.

\para

As we stated above, given  $(a_1,a_2,a_3)\in{\cal C}$, it holds that
$f_3(a_1,a_2,a_3)=h_1(a_1,a_2)-h_2(a_1,a_2)a_3=0$. Thus, in particular,   $(z,r_2(z),r_3(z))\in B\subset {\cal C}$ verifies that $f_3(z,r_2(z),r_3(z))=0$. Hence,
$F_3(z,r_2(z),r_3(z),1)=0$, where $F_3(x_1, x_2, x_3, x_4)$ is the homogeneous polynomial of
$f_3(x_1, x_2, x_3)$.

\para

Taking into account the results in  \cite{paper1}, we have that $r_2(z)=z\varphi_2(1/z)$, where
$\varphi_2(z)\in{{\Bbb C}\ll z\gg}$. Now,  we look for
$\varphi_3(z)\in{{\Bbb C}\ll z\gg}$ such that $r_3(z)=z\varphi_3(1/z)$. This
series must verify that (see statement above)
$$F_3(z,z\varphi_2(1/z),z\varphi_3(1/z),1)=0\qquad \mbox{for $|z|>M$}.$$
We set $z=t^{-1}$, and we get that
$F_3(t^{-1},t^{-1} \varphi_2(t), t^{-1}\varphi_3(t),1)=0$
or  equivalently
$$F_3(1,\varphi_2(t),\varphi_3(t),t)=0.\qquad \qquad \mbox{(I)}$$

\noindent Note that equality (I)  holds for   $|t|<1/M$. That is, equality (I) must be satisfied in a
neighborhood of the infinity point $(1,\varphi_2(0),\varphi_3(0),0)$.

\para

\noindent
At this point, we observe that $F_3$ has the form
$$F_3(x_1,x_2,x_3,x_4)=x_4^{n_1}H_1(x_1,x_2,x_4)-x_4^{n_2}H_2(x_1,x_2,x_4)x_3$$
where $H_i(x_1,x_2,x_4)$ is the homogeneous polynomial of $h_i(x_1,x_2),\,i=1,2$,
and $n_1,n_2\in\mathbb{N}$. Then, we have that
$$F_3(1,\varphi_2(t),\varphi_3(t),t)=t^{n_1}H_1(1,\varphi_2(t),t)-t^{n_2}H_2(1,\varphi_2(t),t)\varphi_3(t)$$
and since (I) must   hold,  we   obtain that
$$\varphi_3(t)=t^{n_1-n_2}\frac{H_1(1,\varphi_2(t),t)}{H_2(1,\varphi_2(t),t)}.$$
 Obviously, $\varphi_3(t)$ can be expressed as a Puiseux
series since ${{\Bbb C}\ll t\gg}$ is a field. Therefore, we conclude that
$B=\{(z,r_2(z), r_3(z))\in {\Bbb C}^3:\,z\in {\Bbb C},\,|z|>M\}$,
where $r_3(z)=z\varphi_3(1/z)$, is an infinity branch of $\cal C$.
\hfill $\Box$

\para
\para

\noindent
In the following, we illustrate the above theorem with an example.

\begin{example}\label{E-infbranches}
Let $\cal C$ be the irreducible space curve defined over $\Bbb C$ by the polynomials
$$f_1(x_1,x_2,x_3)=-x_2^2-2x_1x_3+2x_2x_3-x_1+3,\,\,\mbox{ and}\,\,\,\,\,
f_2(x_1,x_2,x_3)=x_3+x_1x_2-x_2^2.$$ The projection along the
$x_3$-axis, ${\cal C}^p$  is given by the polynomial
$$f^p(x_1,x_2)=x_2^2+x_1-3-2x_2x_1^2+4x_1x_2^2-2x_2^3$$ (this polynomial can be
obtained by computing $\Resultant_{x_3}(f_1, f_2)$; see \cite{SWP}).

\para

By applying the method described in \cite{paper1}, we compute
the infinity branches of ${\cal C}_p$. We obtain the branch
$B_1^p=\{(z,r_{12}(z)):|z|>M_1^p\}$, where
$$r_{12}(z)=\frac{z^{-1}}{2}-\frac{3z^{-2}}{2}+\frac{z^{-3}}{2}-\frac{23z^{-4}}{8}+\frac{37z^{-5}}{8}-\frac{25z^{-6}}{4}+\cdots,$$
that is associated to the infinity point $P_1=(1:0:0)$, and the
branch $B_2^p=\{(z,r_{22}(z)):|z|>M_2^p\}$, where
$$r_{22}(z)=z+\frac{\sqrt{2}z^{1/2}}{2}+\frac{1}{4}+\frac{9\sqrt{2}z^{-1/2}}{32}-\frac{z^{-1}}{4}-\frac{785\sqrt{2}z^{-3/2}}{1024}+\cdots,$$
that is associated to the infinity point $P_2=(1:1:0)$. Note that
$B_2^p$ has ramification index $2$, so
it has two leaves.

\para

Once we have obtained the infinity branches of the projected curve
${\cal C}^p$, we compute the infinity branches of the space curve
${\cal C}$. We use the lift function $h(x_1,x_2)=-x_1x_2+x_2^2$ to
get the third component of these branches (we apply the results in
\cite{Bajaj} to compute $h$). Thus, the infinity branches of the
space curve are $B_1=\{(z,r_{12}(z),r_{13}(z)):|z|>M_1\}$, where
$$r_{13}(z)=h(z,r_{12}(z))=-\frac{1}{2}-\frac{3z^{-1}}{2}-\frac{z^{-2}}{4}+\frac{11z^{-3}}{8}-\frac{15z^{-4}}{8}
+\frac{15z^{-5}}{8}+\cdots$$ and
$B_2=\{(z,r_{22}(z),r_{23}(z)):|z|>M_2\}$, where
$$r_{23}(z)=h(z,r_{22}(z))=\frac{\sqrt{2}z^{3/2}}{2}+\frac{3z}{4}+\frac{17\sqrt{2}z^{1/2}}{32}+\frac{3}{8}
-\frac{897\sqrt{2}z^{-1/2}}{1024}+\cdots.$$
In Figure \ref{F-ejemplo-branches}, we plot the curve $\cal C$ and
some points of the infinity branches $B_1$ and $B_2$.

\begin{figure}[h]
$$
\begin{array}{cc}
\psfig{figure=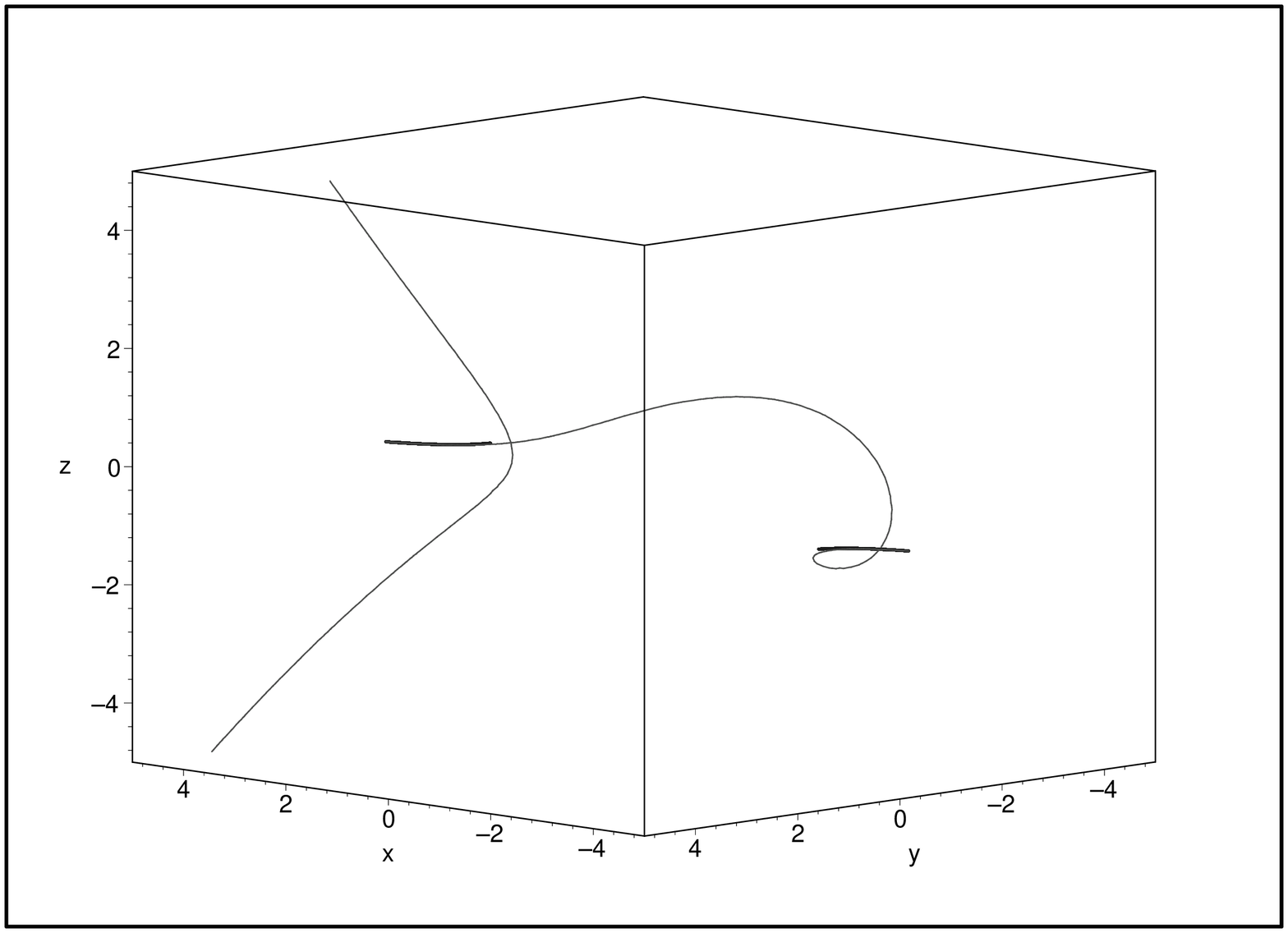,width=5cm,height=5cm,angle=270} &
\psfig{figure=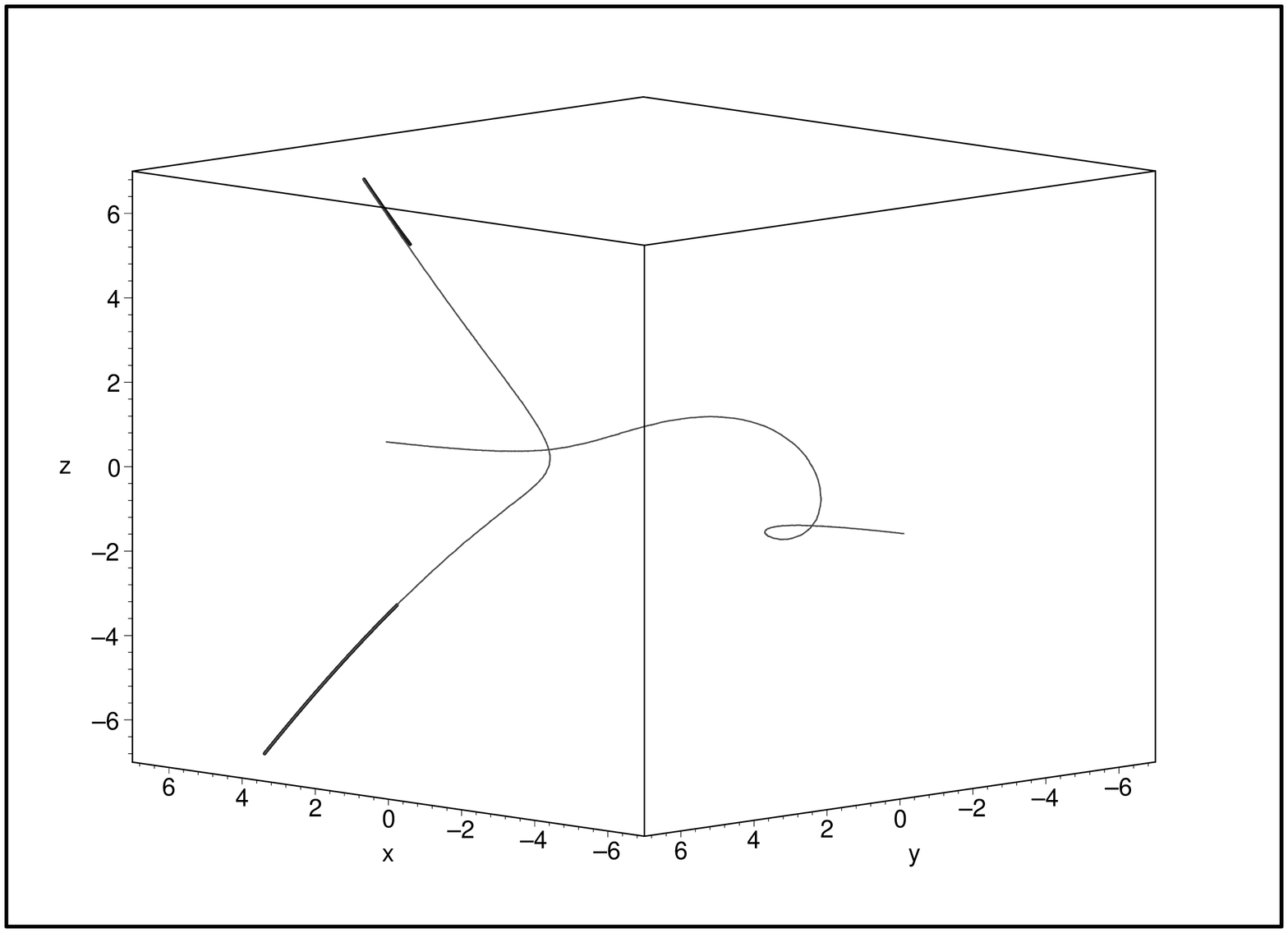,width=5cm,height=5cm,angle=270}
\end{array}
$$ \caption{Curve $\cal C$ and infinity branches $B_1$ (left) and $B_2$
(right).}\label{F-ejemplo-branches}
\vspace*{-0.05cm}
\end{figure}
\end{example}

\section{Computation of an asymptote of a given infinity branch}

In \cite{BlascoPerezII}, we show how some algebraic plane
curves can be approached at infinity by other curves of less degree.
A well-known example is the case of hyperbolas that are curves of degree 2 approached at infinity by two lines (their asymptotes).
Similar situations may also arise when we deal with curves of higher
degree.

\para

\begin{figure}[h]
$$
\begin{array}{cc}
\psfig{figure=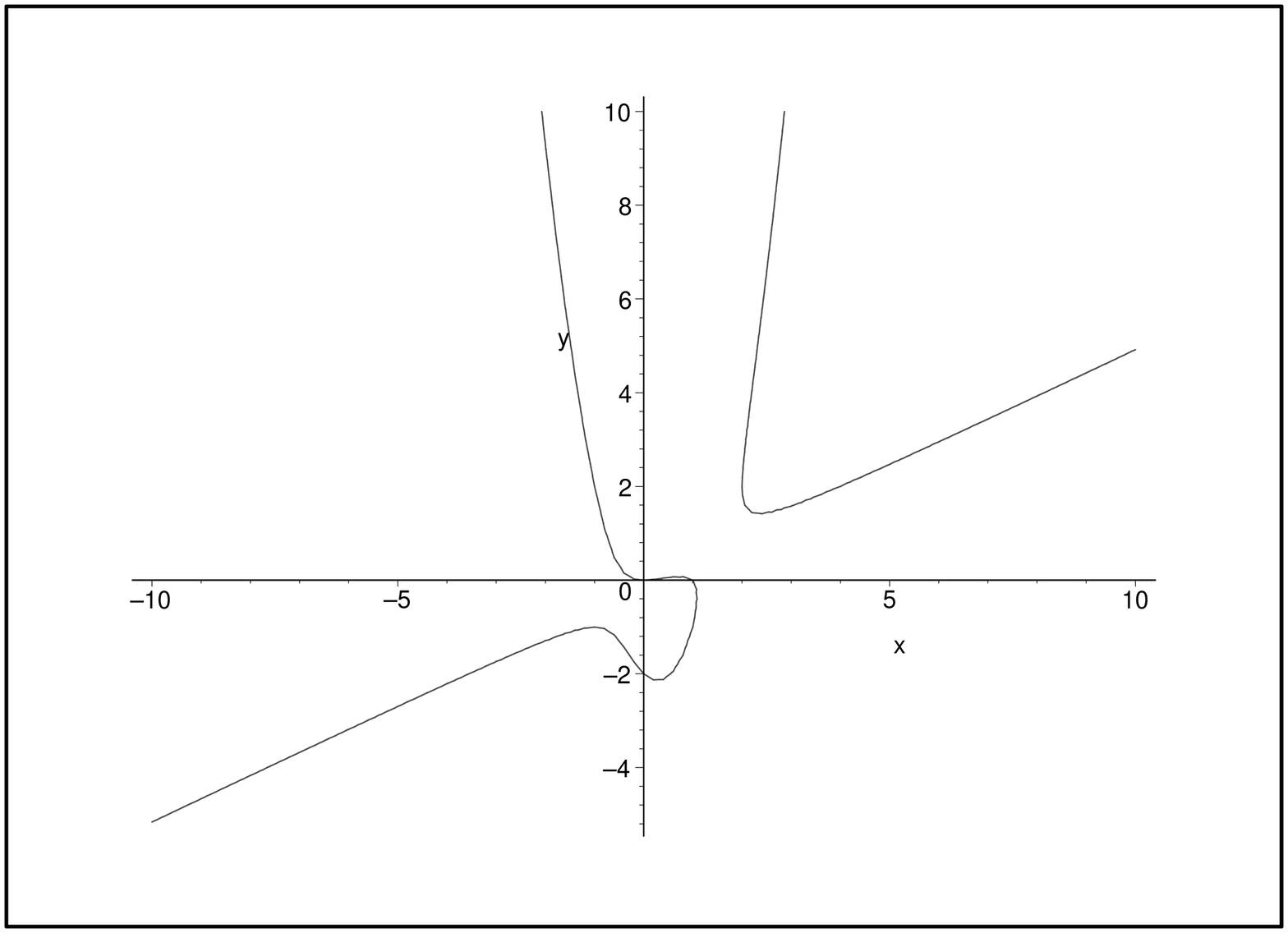,width=4.5cm,height=4.5cm,angle=270} &
\psfig{figure=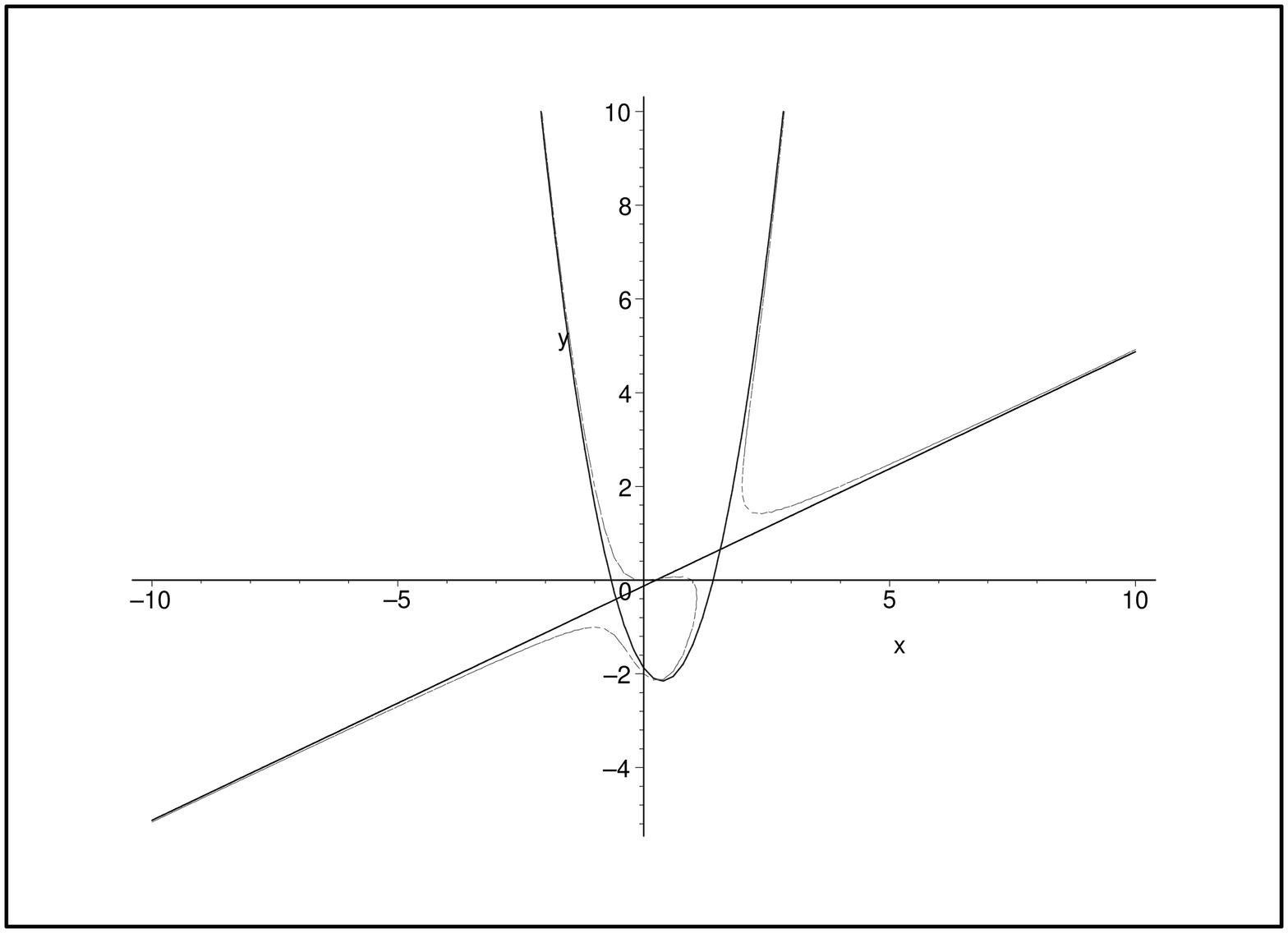,width=4.5cm,height=4.5cm,angle=270}
\end{array}
$$ \caption{Curve $\cal C$ (left)  approached by a parabola and a line (right).}\label{F-ejemplo-asintotasplanas}
\end{figure}

For instance, let $\cal C$ be the plane curve defined by the
equation $-yx-y^2-x^3+2x^2y+x^2-2y=0$. The curve $\cal C$ has degree 3 but it
can be approached at infinity by the parabola $y-2x^2+3/2x+15/8=0$
(see Figure \ref{F-ejemplo-asintotasplanas}). This example leads  us to
introduce the notions of {\em perfect curve} and  {\em
g-asymptote}. Some important properties on these concepts are presented for a given plane curve in Sections 3 and 4 in \cite{BlascoPerezII}. Most of these results can be easily generalized for a given algebraic space curve.


\begin{definition}\label{D-perfect-curve}
A curve of degree $d$ is a {\em perfect curve} if it cannot be approached by
any curve of degree less than $d$.
\end{definition}

\para

A curve that is not perfect can be approached by other curves of
less degree. If these curves are perfect, we call them
{\it g-asymptotes}. More precisely, we have the following definition.

\para

\begin{definition}\label{D-asymptote}
Let ${\cal C}$ be a curve with an infinity branch $B$. A {\em
g-asymptote} (generalized asymptote) of ${\cal C}$ at $B$ is a
perfect curve that approaches ${\cal C}$ at $B$.
\end{definition}

The notion of {\em
g-asymptote} is similar to the classical concept of asymptote. The
difference is that a g-asymptote does not have to be a line, but a
perfect curve. Actually, it is a generalization, since every line is
a perfect curve (this remark follows from Definition
\ref{D-perfect-curve}). Throughout the paper we  refer to
{\em g-asymptote} simply as {\it asymptote}.

\para

\begin{remark}\label{R-minimal-degree}
The degree of an asymptote is less or equal than the
degree of the curve it approaches. In fact, an asymptote of a curve
$\cal C$ at a branch $B$ has minimal degree among all the curves
that approach $\cal C$ at $B$ (see Remark 3 in \cite{BlascoPerezII}).
\end{remark}

\para


 In the following, we prove that every infinity branch
of a given algebraic space curve has, at least, one asymptote and we
show how to obtain it (see Theorem \ref{T-constr-asintota}). Most of
the results introduced bellow to the space case generalize the results presented in
\cite{BlascoPerezII} for the plane case.

\para

Let $\cal C$ be an irreducible space curve implicitly defined by the polynomials
$f_1(x_1,x_2,x_3),\,f_2(x_1,x_2,x_3)\in {\Bbb R}[x_1,x_2,x_3]$, and let $B=\{(z,r_2(z),
r_3(z))\in {\Bbb C}^3:\,z\in {\Bbb C},\,|z|>M\}$ be an infinity
branch of $\cal C$ associated to the infinity point
$P=(1:m_2:m_3:0)$. We know that $r_2$ and $r_3$ are given as
$$r_{2}(z)=m_2z+a_{1,2}z^{-N_{1,2}/N_2+1}+a_{2,2}z^{-N_{2,2}/N_2+1}+a_{3,2}z^{-N_{3,2}/N_2+1} + \cdots$$
$$r_{3}(z)=m_3z+a_{1,3}z^{-N_{1,3}/N_3+1}+a_{2,3}z^{-N_{2,3}/N_3+1}+a_{3,3}z^{-N_{3,3}/N_3+1} + \cdots$$
where $a_{i,2}\not=0$,\, $N_2, N_{i,2}\in {\Bbb
N},\,\,i\geq 1$, \, $0<N_{1,2}<N_{2,2}<\cdots$, and
$a_{i,3}\not=0$,\,  $N_3, N_{i,3}\in {\Bbb N},\,\,i\geq 1$, and
$0<N_{1,3}<N_{2,3}<\cdots$.   Let $N:=\lcm(N_2,N_3)$, and  note that   $\nu(B)=N$.

\para

\begin{lemma}\label{L-degC-N} It holds that
$\deg({\cal C})\geq N.$
\end{lemma}

\noindent\textbf{Proof:} In Section \ref{S-notation}, we show that
there exist $N:=\lcm(N_2,N_3)$ conjugate tuples, $(\varphi_{2}(z), \varphi_{3}(z))$,
which are solutions of the system $g_i(x_2,x_3,x_4)=0$, $i=1,2$. Hence,  the tuples $(z,r_{j,2}(z),r_{j,3}(z))$
with $r_{j,2}(z)=z\varphi_{j,2}(z^{-1})$ and
$r_{j,3}(z)=z\varphi_{j,3}(z^{-1})$ for $j=1,\ldots,N$, are solutions of the system
$f_i(x_1,x_2,x_3)=0$, $i=1,2$. That is, they are points of the curve
$\cal C$.

\para

Then, given $z_0$ such that $|z_0|>M$, we have $N$ intersections
between the curve $\cal C$ and the plane defined by the equation $x_1-z_0=0$ (these points are $(z_0,r_{j,2}(z_0),r_{j,3}(z_0)),\,$  $j=1,\ldots,N$). Thus,
by definition of degree of a space curve (see e.g. \cite{Bajaj91} or \cite{Farouki}), we get that $\deg({\cal C})\geq N$.\hfill $\Box$

\para

\para

\noindent
In the following, we write
{\small \begin{equation}\label{Eq-inf-branchn}
\begin{array}{l}
r_{2}(z)=m_2z+a_{1,2}z^{-\frac{n_{1,2}}{n_2}+1}+\cdots
+a_{\ell_2, 2}z^{-\frac{n_{\ell_2, 2}}{n_2}+1}+a_{\ell_2+1,
2}z^{-\frac{N_{\ell_2+1, 2}}{N_2}+1}+\cdots\\
r_{3}(z)=m_3z+a_{1,3}z^{-\frac{n_{1,3}}{n_3}+1}+\cdots
+a_{\ell_3, 3}z^{-\frac{n_{\ell_3, 3}}{n_3}+1}+a_{\ell_3+1,
3}z^{-\frac{N_{\ell_3+1, 3}}{N_3}+1}+\cdots
\end{array}
\end{equation} }
where $\gcd(n_k,n_{1,k},\ldots,n_{\ell_k,k})=1$, $k=1,2$. That is, we  have simplified the non negative exponents such that
$\gcd(n_k,n_{1,k},\ldots,n_{\ell_k,k})=1$, $k=1,2$. Note that
$0<n_{1,k}<n_{2,k}<\cdots$,  $n_{\ell_k,k}\leq n_k$, and
$N_k<N_{\ell_k+1,k}$. That is,  the terms $a_{j,k}z^{-N_{j,k}/N_k+1}$
with $j\geq \ell_k+1$ have negative exponent.


\para

Under these conditions, we introduce the definition of degree of a branch $B$ as follows:

\para

\begin{definition}\label{D-degreebranch} Let $B=\{(z,r_2(z), r_3(z))\in
{\Bbb C}^3: \,z\in {\Bbb C},\,|z|>M\}$ defined by  (\ref{Eq-inf-branchn})  an
infinity branch associated to   $P=(1:m_2:m_3:0),\,m_j\in {\Bbb C},\,j=1,2$. We say that $n:=\lcm(n_2, n_3)$ is the degree of $B$, and we denote it by $\deg(B)$.
\end{definition}

\begin{remark}\label{R-degC-n} Note that $n_i\leq N_i,\,i=1,2$. Thus, $n=\lcm(n_2, n_3)=\deg(B)\leq N=\lcm(N_2,
N_3)$, and from Lemma \ref{L-degC-N} we get that $\deg({\cal
C})\geq \deg(B).$
\end{remark}


\begin{proposition}\label{P-case3}
Let $\overline{{\cal C}}$ be a curve that approaches ${\cal C}$ at
its infinity branch $B$. It holds that $\deg(\overline{{\cal C}})
\geq \deg(B)$.
\end{proposition}
\noindent\textbf{Proof:} From Theorem \ref{T-curvas-aprox}, we get
that  $\overline{{\cal C}}$ has an infinity branch
$\overline{B}=\{(z,\overline{r}_2(z),\overline{r}_3(z))\in {\Bbb
C}^3:\,z\in {\Bbb C},\,|z|>\overline{M}\}$ convergent with the branch
$B=\{(z,r_2(z),r_3(z))\in {\Bbb C}^3:\,z\in {\Bbb C},\,|z|>M\}$.
From Lemma \ref{L-DistVertical}, we deduce that the terms with non
negative exponent in the series $r_i(z)$ and $\overline{r}_i(z)$,
for $i=2,3$,  are the same, and hence $\overline{B}$ is a branch of
degree $n$ of the form given in (\ref{Eq-inf-branchn}).  Now, the
result follows taking into account Remark \ref{R-degC-n}.\hfill
$\Box$


\subsection{Construction of
asymptotes}\label{S-construction-asymptote}

Let ${\cal C}$ be a space curve with an infinity branch
 $B=\{(z,r_2(z),r_3(z))\in {\Bbb
C}^3:\,z\in {\Bbb C},\,|z|>M\}$. Taking into account the results
presented above, we have  that any curve $\overline{{\cal C}}$
approaching ${\cal C}$ at $B$ has an infinity branch
 $\overline{B}=\{(z,\overline{r}_2(z),\overline{r}_3(z))\in {\Bbb
C}^3:\,z\in {\Bbb C},\,|z|>\overline{M}\}$ such that the terms with
non negative exponent in $r_i(z)$ and $\overline{r}_i(z)$ (for
$i=2,3$) are the same. We consider the series $\tilde{r}_{2}(z)$ and
$\tilde{r}_{3}(z)$, obtained from $r_{2}(z)$ and $r_{3}(z)$ by
removing the terms with negative exponent  (see equation
(\ref{Eq-inf-branchn})). Then, we have that
\begin{equation}\label{Eq-inf-branch3}
\begin{array}{l}
\tilde{r}_{2}(z)=m_2z+a_{1,2}z^{-n_{1,2}/n_2+1}+\cdots +a_{\ell_2,
2}z^{-n_{\ell_2, 2}/n_2+1}\\
\tilde{r}_{3}(z)=m_3z+a_{1,3}z^{-n_{1,3}/n_3+1}+\cdots +a_{\ell_3,
3}z^{-n_{\ell_3, 3}/n_3+1}
\end{array}
\end{equation}
where $a_{j,k},\ldots\in\mathbb{C}\setminus \{0\}$,\,$m_k\in {\Bbb
C}$, $n_k, n_{j,k}\ldots\in\mathbb{N}$,
$\gcd(n_k,n_{1,k},\ldots,n_{\ell,k})=1$, and
$0<n_{1,k}<n_{2,k}<\cdots$. That is, $\tilde{r}_k$ has the same
terms with non negative exponent that $r_k$, and $\tilde{r}_k$ does
not have terms with negative exponent.\\

Let $\widetilde{{\cal C}}$ be the space curve containing the branch
$\widetilde{B}=\{(z,\tilde{r}_2(z),\tilde{r}_3(z))\in {\Bbb
C}^3:\,z\in {\Bbb C},\,|z|>\widetilde{M}\}$. Observe that
\[\widetilde{{\cal Q}}(t)=(t^n, m_2t^n+a_{1,2}t^{r_2(n_2-n_{1,2})}
+\cdots +a_{\ell_2,2}t^{r_2(n_2-n_{\ell_2,2})},\]
\begin{equation}\label{Eq-parametric-case1}
m_3t^n+a_{1,3}t^{r_1(n_3-n_{1,3})} +\cdots
+a_{\ell_3,3}t^{r_3(n_3-n_{\ell_3,3})})\in {\Bbb
C}[t]^3,\end{equation} where $n=\lcm(n_2, n_3)$,\, $r_k=n/n_k$,\, $n_k,
n_{1,k},\ldots,n_{\ell_k,k}\in\mathbb{N}$,\,$0<n_{1,k}<\cdots
n_{\ell_k, k}$ and $\gcd(n_k,n_{1,k},\ldots,n_{\ell_k, k})=1,\,k=2,3$, is a polynomial  parametrization of
$\widetilde{{\cal C}}$. In addition, in Lemma \ref{L-proper-param}, we prove that $\widetilde{{\cal Q}}$ is proper (i.e. invertible).
\para

\begin{lemma}\label{L-proper-param} The parametrization  $\widetilde{{\cal Q}}$ given in (\ref{Eq-parametric-case1}) is proper.
\end{lemma}
\noindent\textbf{Proof:} Let us assume that  $\widetilde{{\cal Q}}$
is not proper. Then,
 there exists $R(t)\in {\Bbb C}[t]$, with $\deg(R)=r>1$,
and ${\cal Q}(t)=(q_1(t),q_2(t),q_3(t))\in {\Bbb C}[t]^3$, such that
${\cal Q}(R)=\widetilde{{\cal Q}}$ (see   \cite{Manocha}). In particular, we get that $q_1(R(t))=t^n$,
which implies that
$$q_1(t)=(t-R(0))^k,\quad \mbox{and}\quad
R(t)=t^r+R(0),\quad rk=n.$$ Let us consider
$R^{\star}(t)=R(t)-R(0)=t^r\in {\Bbb C}[t]$, and
\[{\cal Q}^{\star}(t)={\cal
Q}(t+R(0))=(t^k,{q^{\star}_2}(t),{q^{\star}_3}(t))=\] \[=(t^k,
c_0+c_1t+c_2t^2+\ldots+c_ut^u, d_0+d_1t+d_2t^2+\ldots+d_vt^v) \in
{\Bbb C}[t]^3.\] Then, ${\cal Q}^{\star}(R^{\star})={\cal
Q}(R)=\widetilde{{\cal Q}}$ and, in particular,
\[{q^{\star}_2}(R^\star)={q^{\star}_2}(t^r)=m_2t^n+a_{1,2}t^{r_2(n_2-n_{1,2})}
+\cdots +a_{\ell_2,2}t^{r_2(n_2-n_{\ell_2,2})}\]


\noindent That is,
  $$c_0+c_1t^r+c_2t^{2r}+\ldots+c_ut^{ur}=m_2t^n+a_{1,2}t^{r_2(n_2-n_{1,2})}
+\cdots +a_{\ell_2,2}t^{r_2(n_2-n_{\ell_2,2})}.$$
  From this equality, and taking into account that $r_2=n/n_2=rk/n_2$, we deduce that
  $k/n_2(n_2-n_{i,2})\in\mathbb{Z}$, and thus $kn_{i,2}/n_2\in\mathbb{Z}$ for
  $i=1,\ldots,\ell_2$. This implies that $n_2$ divides $k$ since,
  otherwise, $n_2$ should divide $n_{i,2}$ for
  $i=1,\ldots,\ell_2$, which contradicts the
assumption that $\gcd(n_2,n_{1,2},\ldots,n_{\ell_2, 2})=1$ (see equation
(\ref{Eq-parametric-case1})).

\para

\noindent On the other hand, reasoning similarly  with the
third component, we have that
${q^{\star}_3}(R^\star)={q^{\star}_3}(t^r)=\tilde{q}_3(t)$ and we
get that $n_3$ also divides $k$. Therefore, $k$ is a
common multiple of $n_2$ and $n_3$, which is impossible since $k<n$ (note that $rk=n,\,r>1$)
and $n=\lcm(n_2, n_3)$.\hfill $\Box$ 

\para

From  Lemma \ref{L-proper-param} and using the definition of degree for an implicitly algebraic space
curve (see e.g. \cite{Bajaj91} or \cite{Farouki}),  we obtain the
following lemma.

\begin{lemma}\label{L-case1}
Let $\widetilde{{\cal C}}$ be the plane curve containing the infinity
branch given in (\ref{Eq-inf-branch3}). It holds that $\deg(\widetilde{{\cal C}})=\deg(B)$.
\end{lemma}
\noindent\textbf{Proof:} The intersection of $\widetilde{{\cal C}}$ with a generic plane provides $n$ points since $\widetilde{{\cal C}}$ is parametrized by the proper parametrization  $\widetilde{{\cal Q}}$  that has degree $n$ (see Lemma \ref{L-proper-param}). In addition, we remark that  $n=\deg(B)$ (see Definition
\ref{D-degreebranch}). \hfill $\Box$

\para

In the following theorem, we prove that for any infinity branch $B$ of a  space curve
$\cal C$, there always exists an asymptote that approaches $\cal C$
at $B$. Furthermore, we provide a method to obtain it (see algorithm
{\sf Space Asymptotes Construction}). The proof of this theorem is obtained from Lemmas \ref{L-degC-N} and  \ref{L-case1}, and Proposition \ref{P-case3}. This proof is similar  to the proof of Theorem 2 in \cite{BlascoPerezII}, but for the sake of completeness, we include it.

\para

\begin{theorem}\label{T-constr-asintota}
The curve $\widetilde{{\cal C}}$ is an asymptote of $\mathcal{C}$ at $B$.
\end{theorem}

\noindent\textbf{Proof:} From the construction of
$\widetilde{{\cal C}}$, we have that $\widetilde{{\cal C}}$
approaches $\mathcal{C}$ at $B$. Thus, we  need to show
that  $\widetilde{{\cal
C}}$ cannot be approached by any curve with degree less than
$\deg(\widetilde{{\cal C}})$ (that is, $\widetilde{{\cal C}}$ is perfect).\\
For this purpose, we first note that $\widetilde{{\cal C}}$ has a
polynomial parametrization given by the form in  (\ref{Eq-parametric-case1}). Hence, the unique infinity branch of $\widetilde{{\cal C}}$ is $\widetilde{B}$ (see \cite{Manocha}). In addition, we observe that by construction, $\widetilde{B}$ and $B$ are convergent.\\
Under these conditions, we consider  a plane curve, $\overline{\cal C}$, that approaches $\widetilde{{\cal C}}$
at $\widetilde{B}$. Then,  $\overline{\cal C}$
  has an infinity branch $\overline{B}$ convergent with $\widetilde{B}$ (see Theorem \ref{T-curvas-aprox}). Since   $\widetilde{B}$ and
$B$ are convergent, we deduce that $\overline{B}$ and $B$ are  convergent (see Corollary \ref{C-approaching-curves}) which implies  that $\overline{\cal C}$ approaches $\cal C$ at $B$. Finally, from
 Proposition \ref{P-case3} and Lemma  \ref{L-case1}, we deduce that
$\deg(\overline{\cal C})\geq \deg(\widetilde{{\cal C}})$ and thus, we conclude that $\widetilde{{\cal C}}$ is perfect.\hfill $\Box$

\para

 From
these results, in the following we present an algorithm that computes an
asymptote for each infinity branch of a given space curve.

 \para

 We assume that we have  prepared the input curve $\cal C$,  by means  of a suitable linear change of coordinates if necessary, such that  $(0:a:b:0)$
 ($a\not=0$ or $b\not=0$) is not an infinity point of $\cal C$ (see Remark \ref{R-infinitypoint}). In addition,  we assume that there exists a
birational correspondence between the points of ${\cal C}^p$ and the
points of ${\cal C}$, where ${\cal C}^p$ is the plane curve obtained  by projecting ${\cal C}$ along  the $x_3$-axis (see Section 3).

\para

\begin{center}
\fbox{\hspace*{2 mm}\parbox{13.3cm}{ \vspace*{2 mm} {\bf Algorithm
{\sf Space Asymptotes Construction.}} \vspace*{0.2cm}

\noindent {\sf Given} an irreducible real algebraic space curve $\cal C$
implicitly defined
 by two polynomials $f_1(x_1,x_2,x_3),f_2(x_1,x_2,x_3)\in {\Bbb R}[x_1,x_2,x_3]$,  the algorithm {\sf outputs} an
asymptote for each of its infinity branches.

\begin{itemize}

\item[1.] Compute the projection of $\cal C$ along
the $x_3$-axis. Let  ${\cal C}_p$ be this projection and $f^p(x_1,x_2)$ the implicit polynomial defining ${\cal C}_p$.

\item[2.] Determine the lift function $h(x_1,x_2)$ (see
\cite{Bajaj}).

\item[3.] Compute the  infinity branches of ${\cal C}_p$
by applying Algorithm {\sf Asymptotes Construction} in
\cite{BlascoPerezII}.

\item[4.] For each branch $B^p_i=\{(z,r_{i,2}(z))\in {\Bbb
C}^2:\,z\in {\Bbb C},\,|z|>M^p_{i,2}\}$, $i=1,\ldots, s,$ do:
\begin{itemize}
\item[4.1.] Compute the corresponding infinity branch of $\cal C$:
$$B_i=\{(z,r_{i,2}(z),r_{i,3}(z))\in {\Bbb
C}^3:\,z\in {\Bbb C},\,|z|>M_i\}$$ where
$r_{i,3}(z)=h(z,r_{i,2}(z))$ is given as a Puiseux series.

\item[4.2.] Consider the series $\tilde{r}_{i,2}(z)$ and
$\tilde{r}_{i,3}(z)$  obtained by eliminating the terms with
negative exponent in $r_{i,2}(z)$ and $r_{i,3}(z)$, respectively.
Note that, for $j=2,3$, the series $\tilde{r}_{i,j}$ has the same terms with
non negative exponent that $r_{i,j}$, and $\tilde{r}_{i,j}$ does not
have terms with negative exponent.

\item[4.3.] Return  the asymptote $\widetilde{\cal C}_i$ defined by the proper parametrization (see Lemma \ref{L-proper-param}),
$\widetilde{Q}_i(t)=(t^{n_i},\,\tilde{r}_{i,2}(t^{n_i}),\,\tilde{r}_{i,3}(t^{n_i}))\in
{\Bbb C}[t]^3$, where $n_i=\deg(B_i)$ (see Definition
\ref{D-degreebranch}).
\end{itemize}\end{itemize}
}\hspace{2 mm}}
\end{center}

\para

\begin{remark}\label{R-algoritmo}\begin{itemize}
\item[1.] The implicit polynomial $f^p(x_1,x_2)$  defining ${\cal C}_p$ (see step 1) can be computed as $f^p(x_1,x_2)=\resultant_{x_3}(f_1,f_2)$ (see Section 4.5 in \cite{SWP}).
\item[2.] Since   we have assumed that the given algebraic space curve $\cal C$ only has infinity points of the form $(1 : m_2 : m_3 : 0)$ (see Remark \ref{R-infinitypoint}), we have that $(0:m:0)$ is not an infinity point of the plane curve ${\cal C}_p$ and thus,   Algorithm {\sf Asymptotes Construction} in
\cite{BlascoPerezII} (see step 3) can be applied.
    \item[3.] Reasoning as in the correctness of the algorithm {\sf Asymptotes Construction} in \cite{BlascoPerezII}, one may prove that the algorithm
{\sf  Space Asymptotes Construction} outputs an asymptote $\widetilde{\cal C}$ that is independent of the leaf chosen to define  the branch $B$ (see Section 2).
\end{itemize}
\end{remark}

\para

In the following example, we  illustrate algorithm {\sf Space Asymptotes Construction.} \para

\begin{example}
Let $\cal C$ be the algebraic space curve over $\Bbb C$ introduced in Example \ref{E-infbranches}. The curve $\cal C$ is defined by the polynomials
$$f_1(x_1,x_2,x_3)=-x_2^2-2x_1x_3+2x_2x_3-x_1+3,\,\mbox{ and}\,\,
f_2(x_1,x_2,x_3)=x_3+x_1x_2-x_2^2.$$ In Example \ref{E-infbranches},
we show  that $\cal C$ has two infinity branches given by:
$$B_1=\{(r_{11}(z),r_{12}(z),r_{13}(z)):|z|>M_1\},\qquad \mbox{where}$$
$$\begin{array}{l}r_{11}(z)=z,\\
\\
\displaystyle r_{12}(z)=\frac{z^{-1}}{2}-\frac{3z^{-2}}{2}+\frac{z^{-3}}{2}-\frac{23z^{-4}}{8}+\frac{37z^{-5}}{8}-\frac{25z^{-6}}{4}+\cdots,\\
\\
\displaystyle
r_{13}(z)=-\frac{1}{2}-\frac{3z^{-1}}{2}-\frac{z^{-2}}{4}+\frac{11z^{-3}}{8}-\frac{15z^{-4}}{8}
+\frac{15z^{-5}}{8}+\cdots,
\end{array}$$ and
$$B_2=\{(r_{21}(z),r_{22}(z),r_{23}(z)):|z|>M_2\},\qquad \mbox{where}$$
$$\begin{array}{l}r_{21}(z)=z,\\
\\
\displaystyle
r_{22}(z)=z+\frac{\sqrt{2}z^{1/2}}{2}+\frac{1}{4}+\frac{9\sqrt{2}z^{-1/2}}{32}-\frac{z^{-1}}{4}-\frac{785\sqrt{2}z^{-3/2}}{1024}+\cdots,\\
\\
\displaystyle
r_{23}(z)=\frac{\sqrt{2}z^{3/2}}{2}+\frac{3z}{4}+\frac{17\sqrt{2}z^{1/2}}{32}+\frac{3}{8}
-\frac{897\sqrt{2}z^{-1/2}}{1024}+\cdots.
\end{array}$$

These branches were obtained by applying steps 1, 2, 3, and 4.1
of Algorithm {\sf Space Asymptotes Construction}. Now we apply
step 4.2, and we compute the series $\tilde{r}_{i,j}(z)$ by removing
the terms with negative exponent from the series $r_{i,j}(z)$,
$i=1,2$, $j=1,2,3$. We get:
$$\begin{array}{lll}\tilde{r}_{11}(z)=z, & \hspace{2cm} & \tilde{r}_{21}(z)=z,\\
\\
\displaystyle
\tilde{r}_{12}(z)=0, & \hspace{2cm} & \displaystyle\tilde{r}_{22}(z)=z+\frac{\sqrt{2}z^{1/2}}{2}+\frac{1}{4},\\
\\
\displaystyle \tilde{r}_{13}(z)=-\frac{1}{2}, & \hspace{2cm} &
\displaystyle\tilde{r}_{23}(z)=\frac{\sqrt{2}z^{3/2}}{2}+\frac{3z}{4}+\frac{17\sqrt{2}z^{1/2}}{32}+\frac{3}{8}.
\end{array}$$

\begin{figure}[h]
$$
\begin{array}{cc}
\psfig{figure=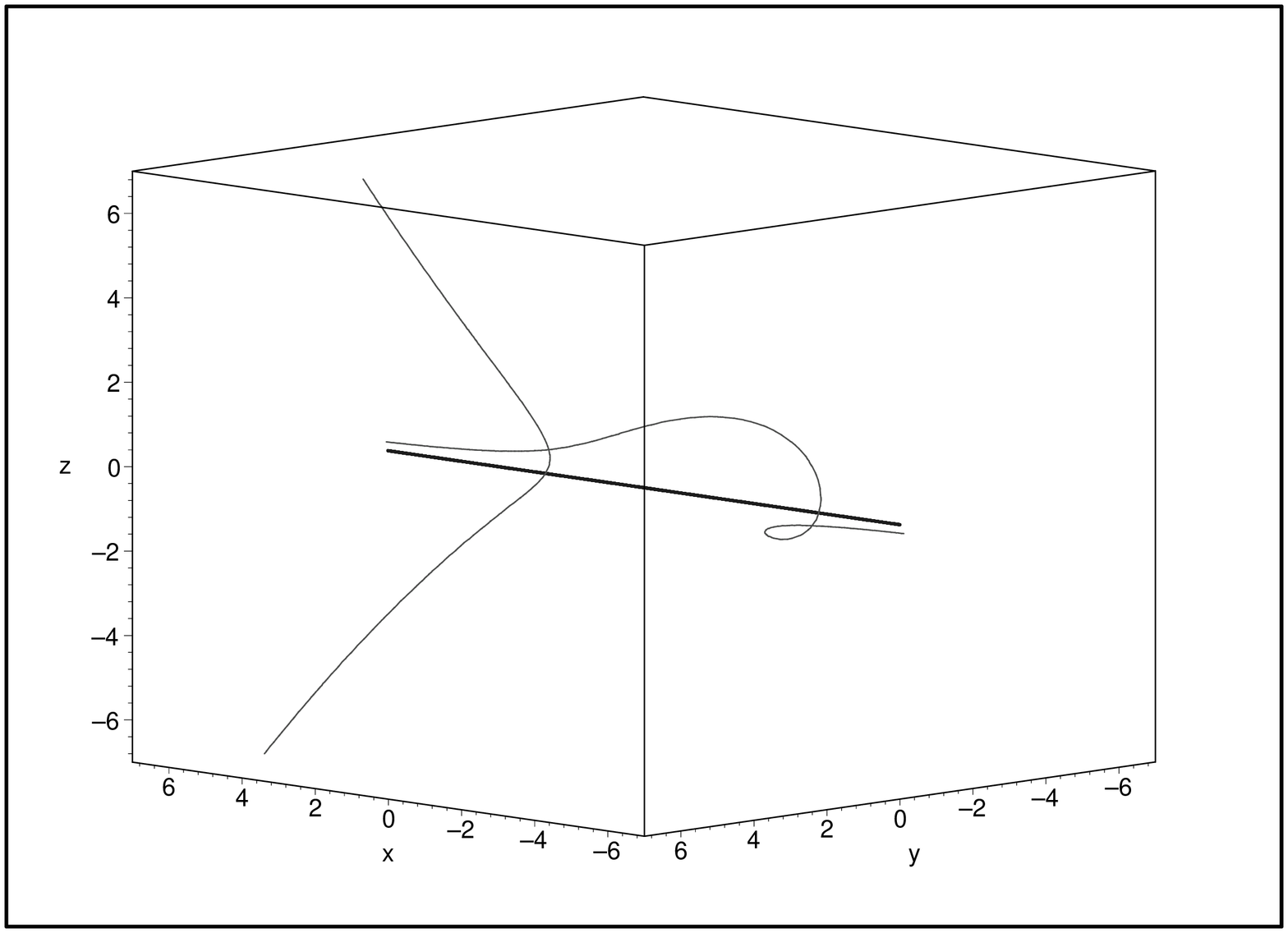,width=5cm,height=5cm,angle=270} &
\psfig{figure=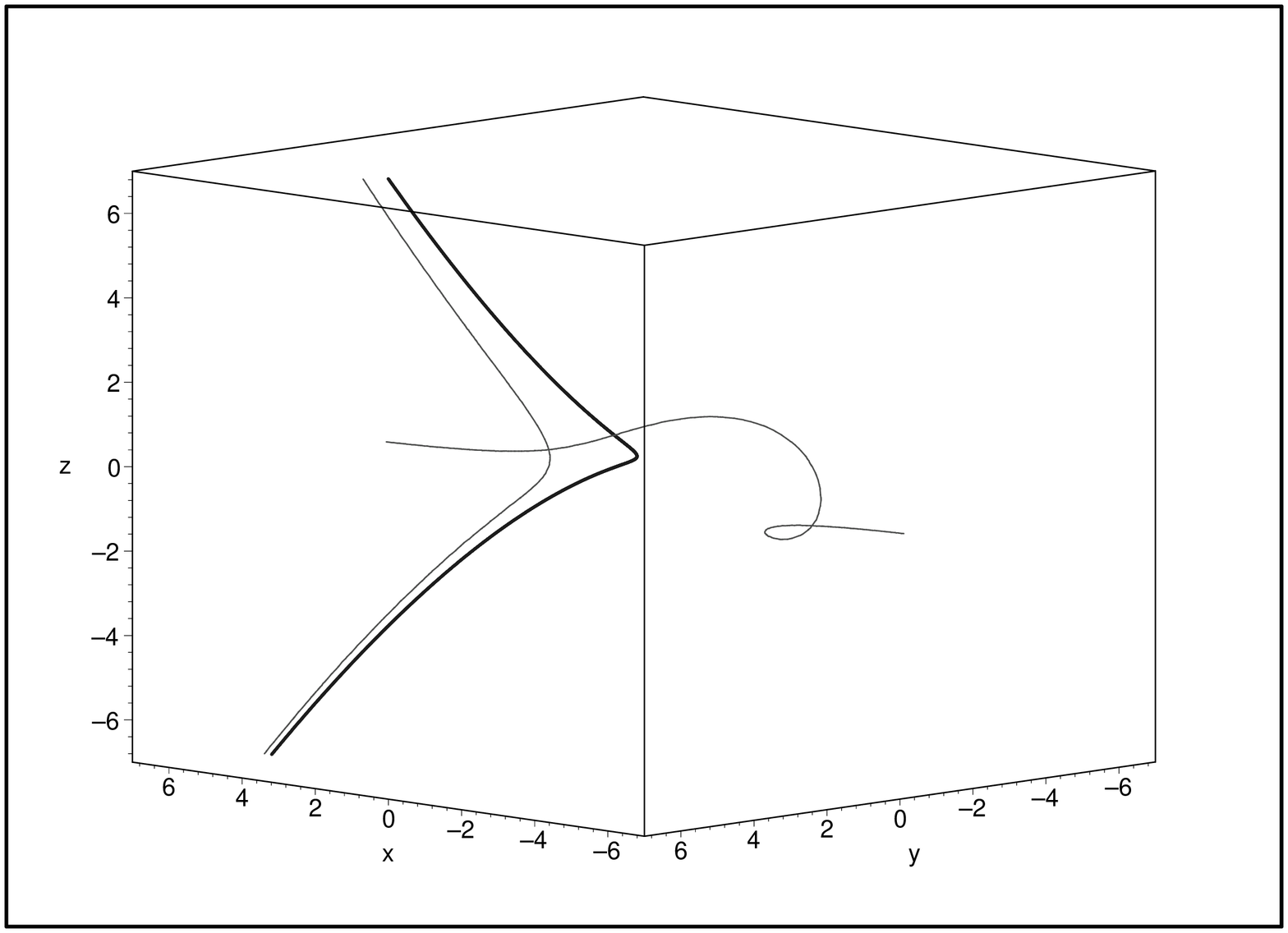,width=5cm,height=5cm,angle=270}
\end{array}
$$ \caption{Curve $\cal C$ approached by asymptotes $\widetilde{{\cal C}}_1$ (left) and $\widetilde{{\cal C}}_2$ (right).}\label{F-ejemplo-asintotas}
\end{figure}

\noindent Thus, in step 4.3, we obtain:

$$\widetilde{Q}_1(t)=(t,\,\tilde{r}_{1,2}(t),\,\tilde{r}_{1,3}(t))=(t,0,-1/2),\qquad \mbox{and}$$
 $$\widetilde{Q}_2(t)=(t^2,\,\tilde{r}_{2,2}(t^2),\,\tilde{r}_{2,3}(t^2))=\left(t^2,t^2+\frac{\sqrt{2}t}{2}+\frac{1}{4},\frac{\sqrt{2}t^3}{2}+\frac{3t^2}{4}+\frac{17\sqrt{2}t}{32}+\frac{3}{8}\right).$$

\noindent $\widetilde{Q}_1$ and $\widetilde{Q}_2$ are proper
parametrizations (see Lemma \ref{L-proper-param}) of the asymptotes $\widetilde{{\cal C}}_1$ and
$\widetilde{{\cal C}}_2$, which approach $\cal C$ at its infinity
branches $B_1$ and $B_2$, respectively.

\para

\noindent
In Figure \ref{F-ejemplo-asintotas}, we plot the curve $\cal C$ and
its asymptotes $\widetilde{{\cal C}}_1$ and $\widetilde{{\cal
C}}_2$.

\end{example}

\section{Asymptotes of a parametric curve}

Throughout this paper, we have dealt with real algebraic space curves defined
implicitly by two polynomials. In this section, we present a method
to compute infinity branches and asymptotes of rational curves from
their parametric representation (without implicitizing).


\para

Thus, in the following, we deal with real space curves  defined parametrically. However, the method described  can be trivially
applied to the case of parametric real plane curves and in general, for a rational parametrization of a curve in the $n$-dimensional space. 
Similarly as in the previous sections,   we  work
over   $\Bbb C$, but we assume that the curve has infinitely many points in the affine plane over $\Bbb R$ and then, the curve has a real parametrization  (see Chapter 7 in \cite{SWP}).

\para

Under these conditions, in the following, we consider a real space curve ${\cal C}$  defined by the parametrization
\[{\cal P}(s)=(p_1(s),p_2(s),p_3(s))\in {\Bbb R}(s)^3\setminus{\Bbb R}^3,\quad p_i(s)=p_{i1}(s)/p(s),\,i=1,2,3.\]
We assume that   we have  prepared the input curve ${\cal C}$, by means
of a suitable linear change of coordinates (if necessary) such that
$(0:a:b:0)$ ($a\not=0$ or $b\not=0$) is not an infinity point  (see Remark \ref{R-infinitypoint}). Note that this implies that  $\deg(p_1)\geq 1$.

\para

Observe  that if ${\cal C}^*$ represents  the projective curve associated to
${\cal C}$, we have that a parametrization of ${\cal C}^*$ is given by
${\cal P}^*(s)=(p_{11}(s):p_{21}(s):p_{31}(s):p(s))$ or, equivalently,
$${\cal P}^*(s)=\left(1:\frac{p_{21}(s)}{p_{11}(s)}:\frac{p_{31}(s)}{p_{11}(s)}:\frac{p(s)}{p_{11}(s)}\right).$$

\para

\noindent
\underline{\sf A method to construct the asymptotes of $\cal C$.}\\


In order to compute the asymptotes of $\cal C$, first we need to determine the  infinity branches of $\cal C$. That is, the sets $B=\{(z:r_{2}(z):r_{3}(z)):\,z\in {\Bbb
C},\,|z|>M\},$ where $r_{j}(z)=z\varphi_{j}(z^{-1}),\, j=2,3$. For this purpose, we note that from Section \ref{S-notation}, we have that $F_i(1:\varphi_{2}(t):\varphi_{3}(t):t)=0$ around $t=0$,  where $F_i,\,i=1,2$ are the polynomials defining implicitly ${\cal C}^*$.  Observe that in this section, we are given the parametrization ${\cal P}^*$ of ${\cal C}^*$ and  then,  $F_i({\cal P}^*(s))=F_i(1:\frac{p_{21}(s)}{p_{11}(s)}:\frac{p_{31}(s)}{p_{11}(s)}:\frac{p(s)}{p_{11}(s)})=0$. Thus, intuitively speaking, in order to compute the  infinity branches of $\cal C$, and in particular the series $\varphi_{j},\,j=2,3$, one needs to rewrite the parametrization ${\cal P}^*(s)=\left(1:\frac{p_{21}(s)}{p_{11}(s)}:\frac{p_{31}(s)}{p_{11}(s)}:\frac{p(s)}{p_{11}(s)}\right)$ in the form $(1:\varphi_{2}(t):\varphi_{3}(t):t)$ around $t=0$. For this purpose, the idea is to look for a value of the parameter $s$, say $\ell(t)\in {\Bbb C}\ll
t\gg$, such that ${\cal P}^*(\ell(t))=(1:\varphi_{2}(t):\varphi_{3}(t):t)$ around $t=0$.\\

Hence, from the above reasoning, we deduce that first, we have to consider the equation $p(s)/p_{11}(s)=t$ (or equivalently,  $p(s)-tp_{11}(s)=0$), and we solve it in the variable $s$ around $t=0$ (note that $\deg(p_1)\geq 1$). From Puiseux's Theorem,   there
exist solutions $\ell_1(t),\ell_2(t),\ldots,\ell_k(t)\in {\Bbb C}\ll
t\gg$ such that,
$p(\ell_i(t))-tp_{11}(\ell_i(t))=0,\,i=1,\ldots,k,$ in a neighborhood of $t=0$.\\

Thus,   for each $i=1,\ldots,k$, there exists  $M_i\in {\Bbb R}^+$ such that
the points $(1:\varphi_{i,2}(t):\varphi_{i,3}(t):t)$ or  equivalently, the
points $(t^{-1}:t^{-1}\varphi_{i,2}(t):t^{-1}\varphi_{i,3}(t):1)$, where
\begin{equation}\label{Eq-psi}\varphi_{i,j}(t)=\frac{p_{j,1}(\ell_i(t))}{p_{11}(\ell_i(t))},\quad
j=2,3,\end{equation} are in ${\cal C}^*$ for $|t|<M_i$ (note that ${\cal P}^*(\ell(t))\in {\cal C}^*$ since ${\cal P}^*$ is a parametrization of ${\cal C}^*$). Observe  that $\varphi_{i,j}(t),\, j=2,3$ are
Puiseux series, since $p_{j,1}(\ell_i(t)),\, j=2,3$ and
$p_{11}(\ell_i(t))$ can be written as Puiseux series and ${\Bbb C}\ll t\gg$ is a
field.\\

Finally, we set $z=t^{-1}$. Then, we have that the points
$(z:r_{i,2}(z):r_{i,3}(z))$, where $r_{i,j}(z)=z\varphi_{i,j}(z^{-1}),
j=2,3$, are in ${\cal C}$ for $|z|>M_i^{-1}$. Hence, the infinity
branches of $\cal C$ are the sets
$$B_i=\{(z:r_{i,2}(z):r_{i,3}(z)):\,z\in {\Bbb
C},\,|z|>M_i^{-1}\},\quad i=1,\ldots,k.$$

\para

\begin{remark}\label{R-calculo-r-directo}
Note that the series $\ell_i(t)$ satisfies that
$p(\ell_i(t))/p_{11}(\ell_i(t))=t$, for $i=1,\ldots,k$. Then, from equality
(\ref{Eq-psi}), we have that for $j=2,3$
$$\varphi_{i,j}(t)=\frac{p_{j,1}(\ell_i(t))}{p(\ell_i(t))}t=p_j(\ell_i(t))t,\quad\mbox{and}\quad r_{i,j}(z)=z\varphi_{i,j}(z^{-1})=p_j(\ell_i(z^{-1})).$$
\end{remark}

\para

Once we have the infinity branches, we can compute an asymptote for
each of them by simply removing the terms with negative exponent
from $r_{i,2}$ and $r_{i,3}$ (see Subsection
\ref{S-construction-asymptote}).

\para

The following algorithm computes the infinity branches of a given
parametric space curve and provides an asymptote for each of them.
We remind that the input curve $\cal C$ is prepared such that  $(0:a:b:0)$ ($a\not=0$ or $b\not=0$) is not an infinity point of
${\cal C}^*$  (see Remark \ref{R-infinitypoint}).


\para

\begin{center}
\fbox{\hspace*{2 mm}\parbox{13.2cm}{ \vspace*{2 mm} {\bf Algorithm
{\sf Space Asymptotes Construction-Parametric Case.}}
\vspace*{0.2cm}

\noindent {\sf Given} a rational irreducible real  algebraic space curve
$\cal C$ defined by a   parametrization
${\cal P}(s)=(p_1(s),p_2(s),p_3(s))\in \mathbb{R}(s)^3,$\,$p_j(s)=p_{j1}(s)/p(s),\,j=1,2,3$,  the algorithm {\sf outputs}  one asymptote
for each of its infinity branches.
\begin{itemize}

\item[1.]  Compute the Puiseux solutions of $p(s)-tp_{11}(s)=0$
around $s=0$. Let them be $\ell_1(t),\ell_2(t),\ldots,\ell_k(t)\in
{\Bbb C}\ll t\gg$.

\item[2.]  For each   $\ell_i(t)\in
{\Bbb C}\ll t\gg$, $i=1,\ldots, k,$ do:
\begin{itemize}

\item[2.1.] Compute the corresponding infinity branch of $\cal C$:
$$B_i=\{(z,r_{i,2}(z),r_{i,3}(z))\in {\Bbb
C}^3:\,z\in {\Bbb C},\,|z|>M_i\},\quad \mbox{where}$$
$r_{i,j}(z)=p_j(\ell_i(z^{-1})), \,j=2,3$ is given as Puiseux
series (see Remark \ref{R-calculo-r-directo}).

\item[2.2.] Consider the series $\tilde{r}_{i,2}(z)$ and
$\tilde{r}_{i,3}(z)$  obtained by eliminating the terms with
negative exponent in $r_{i,2}(z)$ and $r_{i,3}(z)$, respectively.
Note that, for $j=2,3$, the series $\tilde{r}_{i,j}$ has the same terms with
non negative exponent that $r_{i,j}$, and $\tilde{r}_{i,j}$ does not
have terms with negative exponent.

\item[2.3.] Return  the asymptote $\widetilde{\cal C}_i$ defined by the proper parametrization (see Lemma \ref{L-proper-param}),
$\widetilde{Q}_i(t)=(t^{n_i},\,\tilde{r}_{i,2}(t^{n_i}),\,\tilde{r}_{i,3}(t^{n_i}))\in
{\Bbb C}[t]^3$, where $n_i=\deg(B_i)$ (see Definition
\ref{D-degreebranch}).
\end{itemize}\end{itemize}
}\hspace{2 mm}}
\end{center}

\para

\begin{remark} We note that:
\begin{enumerate}
\item  In step $1$ of the algorithm,  some of the  solutions $\ell_1(t),\ell_2(t),\ldots,\ell_k(t)\in
{\Bbb C}\ll t\gg$ might belong to the same conjugation class. Thus, we only consider one
solution for each of these classes.
\item Reasoning as in statement $3$ in Remark \ref{R-algoritmo}, one also gets that the algorithm
{\sf   Space Asymptotes Construction-Parametric Case} outputs an asymptote $\widetilde{\cal C}$ that is independent of the solutions $\ell_1(t),\ell_2(t),\ldots,\ell_k(t)\in
{\Bbb C}\ll t\gg$ chosen in step $1$ (see statement $1$ above), and of the leaf chosen to define  the branch $B$.
\end{enumerate}
\end{remark}

\para

In the following example, we study a parametric space curve with only one
infinity branch. We use algorithm {\sf Space Asymptotes
Construction-Parametric Case} to obtain the branch and compute an
asymptote for it.

\para

\begin{example}
Let ${\cal C}$ be the space curve defined by the
parametrization
$${\cal P}(s)=\left(\frac{-1+s^2}{s^3},\frac{-1+s^2}{s^2},\frac{1}{s}\right)\in {\Bbb R}(s)^3.$$
\noindent{\sf Step 1:} We compute the solutions of the equation
$$p(s)-tp_{11}(s)=s^3-t(-1+s^2)=s^3-ts^2+t=0$$
around $t=0$. There is only one solution that is given by the
Puiseux series (see Proposition \ref{P-param})
$$\ell(t)=(-t)^{1/3}+1/3t+1/9(-t)^{5/3}-2/81(-t)^{7/3}+2/729(-t)^{11/3}+\cdots$$
(note that $\ell(t)$ represents a conjugation class composed by
three conjugated series; one of them
is real and the other two are complex).

\para

\noindent{\sf Step 2:}

\begin{itemize}
\item[]  \mbox{\sf {Step 2.1}:} We compute (see Proposition \ref{P-param})
$$r_2(z)=p_2(\ell(z^{-1}))=-z^{2/3}+1/3-1/9z^{-2/3}+2/81z^{-4/3}-2/729z^{-8/3}+\cdots$$
$$r_3(z)=p_3(\ell(z^{-1}))=-z^{1/3}-1/3z^{-1/3}+1/81z^{-5/3}-1/243z^{-7/3}+\cdots.$$
The curve has only one infinity branch given by
$$B=\{(z,r_2(z),r_3(z)):\,z\in {\Bbb
C},\,|z|>M\}$$
for some $M\in {\Bbb R}^+$ (note that this branch has three leaves; one of them
is real and the other two are complex).

\item[]  \mbox{\sf {Step 2.2}:} We obtain $\tilde{r}_2(z)$ and $\tilde{r}_3(z)$ by
eliminating the terms with negative exponent in $r_{2}(z)$ and
$r_{3}(z)$ respectively:
$$\tilde{r}_2(z)=-z^{2/3}+1/3\quad\text{ and }\quad\tilde{r}_3(z)=-z^{1/3}.$$

\item[]  \mbox{\sf {Step 2.3}:} The input curve $\cal C$ has an asymptote $\widetilde{{\cal C}}$ at
$B$ that can be polynomially parametrized by:
$$\widetilde{Q}(t)=(t^3,\,\tilde{r}_{2}(t^3),\,\tilde{r}_{3}(t^3))=(t^3,-t^2+1/3,-t).$$
\end{itemize}

\noindent In Figure \ref{F-ejemplo-parametrico}, we plot the curve
$\cal C$, the infinity branch $B$, and the asymptote
$\widetilde{{\cal C}}$.

\begin{figure}[h]
$$
\begin{array}{ccc}
\psfig{figure=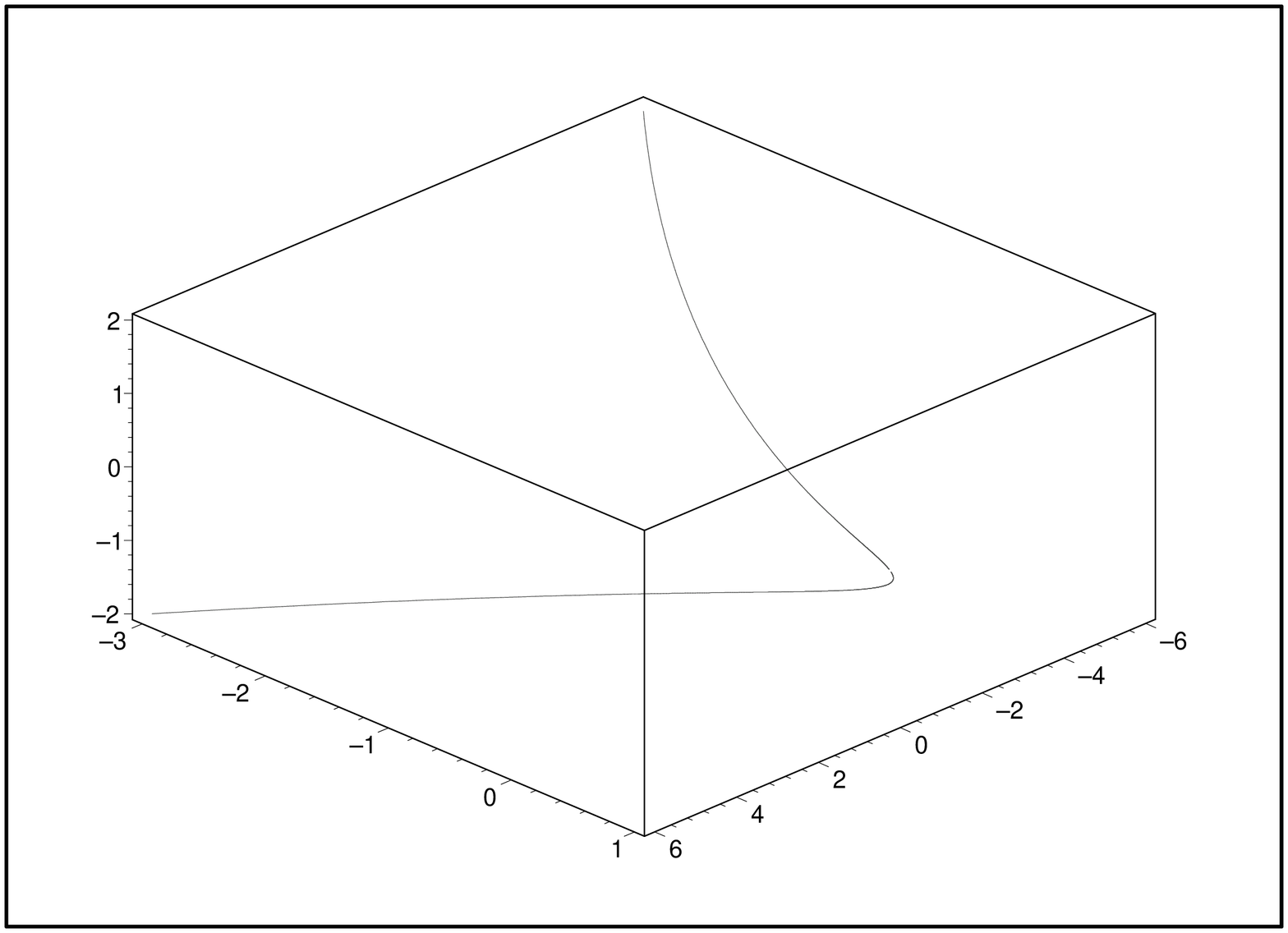,width=4.3cm,height=4.3cm,angle=270} &
\psfig{figure=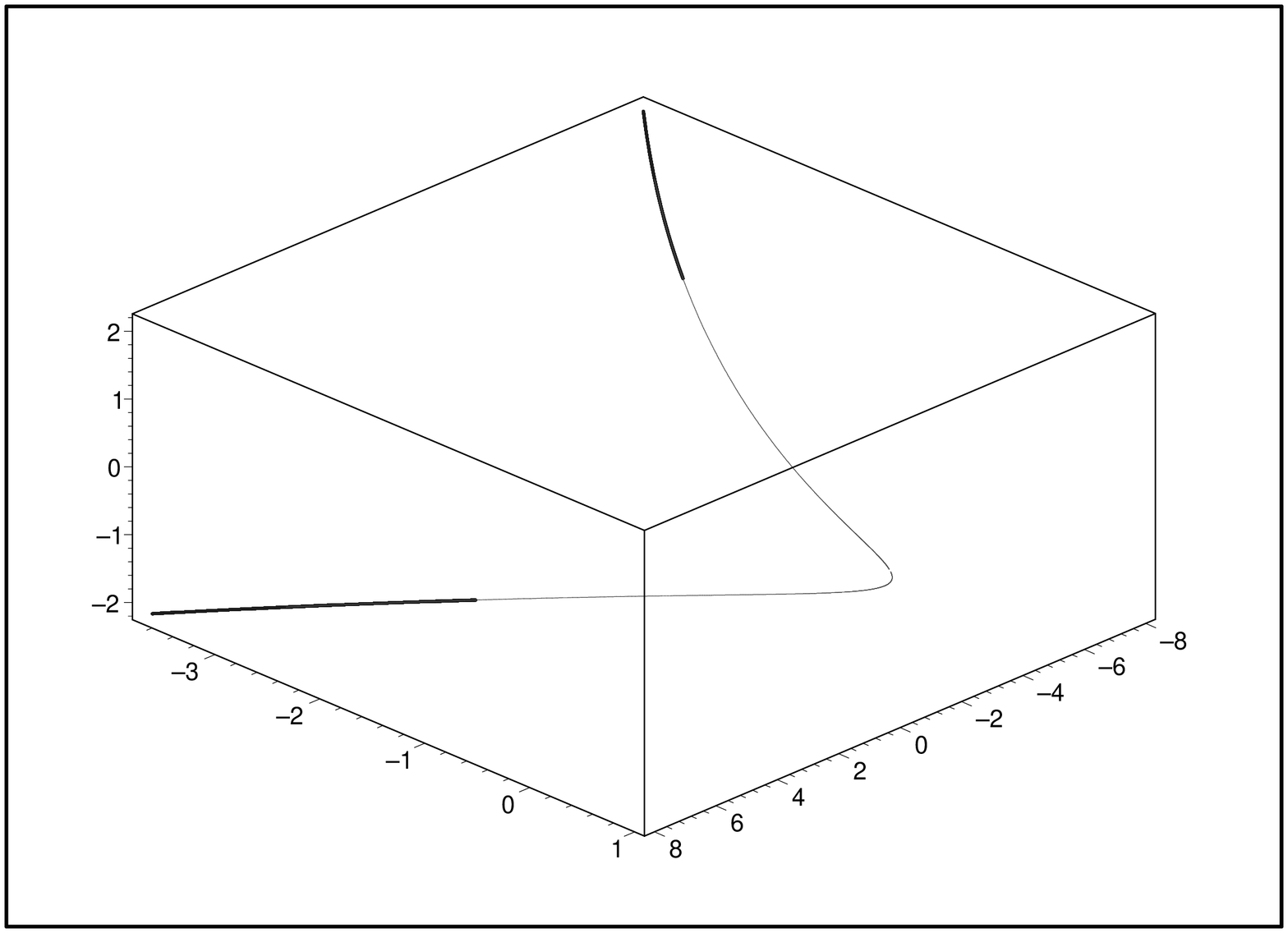,width=4.3cm,height=4.3cm,angle=270} &
\psfig{figure=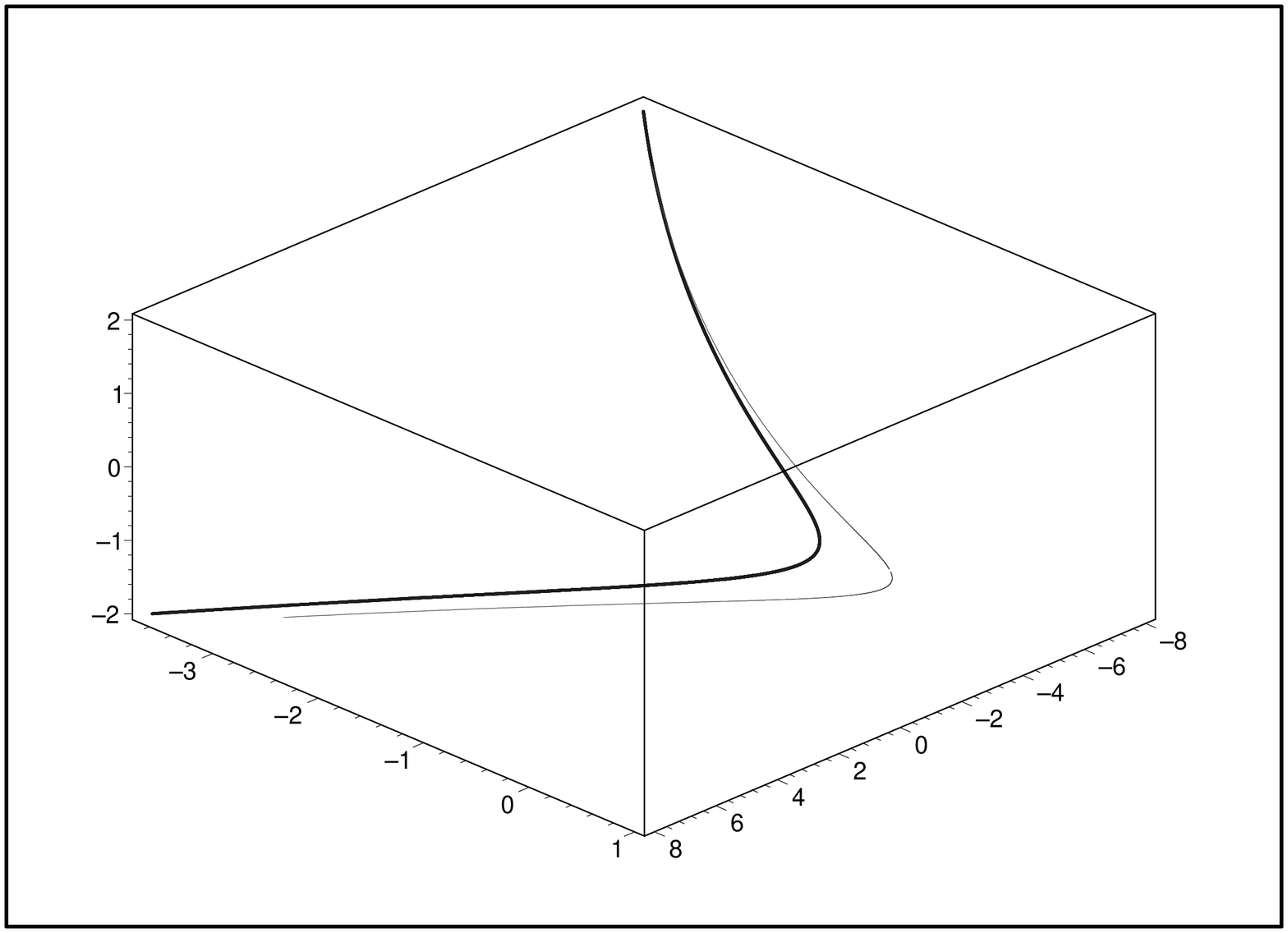,width=4.3cm,height=4.3cm,angle=270}
\end{array}
$$ \caption{Curve $\cal C$ (left), infinity branch
$B$ (center) and asymptote $\widetilde{{\cal C}}$
(right)}\label{F-ejemplo-parametrico}
\end{figure}
\end{example}

\para

\noindent
\underline{\sf Correctness}. \para

The application of the algorithm {\sf Space
Asymptotes Construction-Parametric Case} presents some technical
difficulties since infinite series are
involved. In particular, when we compute the series ${\ell}_i$ in step
1, we cannot handle its infinite terms so it must be truncated,
which may distort the computation of the series $r_{i,j}$ in step 2.
However, this distortion may not affect to all the terms in
$r_{i,j}$. In fact, the number of affected terms depends on the
number of terms considered in ${\ell}_i$. Nevertheless, note that we do not
need to know the full expression of $r_{i,j}$ but only the terms
with non negative exponent. Proposition \ref{P-param}  states that the
terms with non negative exponent in  $r_{i,j}$ can be
obtained from a finite number of terms considered in $\ell_i$. In fact,
it provides a lower bound for the number of terms needed in $\ell_i$.

\para

\begin{proposition} \label{P-param}   Let $\ell(z)\in
{\Bbb C}\ll z\gg$ be a solution obtained in step $1$ of  the
algorithm {\sf Space Asymptotes Construction-Parametric Case}.  Let
$B=\{(z,r_{2}(z),r_{3}(z))\in {\Bbb C}^3:\,z\in {\Bbb
C},\,|z|>M\},\,$ $r_{j}(z)=p_j(\ell(z^{-1})), \,j=2,3,$ be the
infinity branch of $\cal C$  obtained in step $2.1$ of  the
algorithm {\sf Space Asymptotes Construction-Parametric Case}. It
holds that the terms with non negative exponent in  $r_2$ and $r_3$
can be obtained from the computation of $2\deg(p_1)+1$ terms of
$\ell$.
\end{proposition}

\vspace*{2mm}

\noindent {\bf Proof.}  We prove the proposition for $r_2$ (similarly, one gets the result for $r_3$). For this purpose, we
  write $\ell(z)$ as  \[\ell(z):=b_0 +b_{1}z^{-1/N}+\cdots
+b_{k}z^{-{k}/N}+B(z),\quad B(z)=\sum_{j=1}^\infty
a_{j}z^{j/N},\,\quad N \in\mathbb{N}^+,\] $a_i, b_{i}\in\mathbb{C}$, and we consider  ${\ell}^*(z):=\ell(z^N)= \nu/z^k$ where
\[\nu:=b_0 z^{k}+b_{1}z^{k-1}+\cdots
+b_{k-1}z+b_{k}+z^{k}B(z^N),\qquad B(z^N)=\sum_{j=1}^\infty
a_{j}z^{j}.\]

 Note that the terms with non negative exponent in $r_2(z)$
are the terms with non positive exponent in $r_2(1/z)$. In addition,
these terms are the terms with non positive exponent in
$r_2(1/z^N)$. On the other hand, $r_{2}(z)=p_2(\ell(z^{-1}))$ so
$r_2(1/z^N)=p_2({\ell}^*(z))$. Therefore, we need to determine the
terms with non
positive exponent in $p_2({\ell}^*(z))$.\\

\noindent Now, we distinguish two different cases:
\begin{enumerate}
 \item  Let us assume that $\ell(z)$ has terms with negative exponent and thus,  we assume w.l.o.g. that $b_k\not=0,\,k>0$.  Thus,
 \[p_2({\ell}^*(z))=\frac{p_{2,1}(\nu/z^k)}{p(\nu/z^k)}=\frac{\bar{p}_{2,1}(z)}{z^{k(m-n)}\bar{p}(z)},\qquad m:=\deg(p_{2,1}),\,n:=\deg(p),\]
\[\bar{p}_{2,1}(z)=c_m\nu^{m}+c_{m-1}z^{k}\nu^{m-1}+c_{m-2}z^{2k}\nu^{m-2}+\cdots+c_{0}z^{km},\,c_m\not=0\]
\[\bar{p}(z)=d_n\nu^{n}+d_{n-1}z^{k}\nu^{n-1}+d_{n-2}z^{2k}\nu^{n-2}+\cdots+d_{0}z^{kn},\,d_n\not=0.\]
Under these conditions,  the generalized series
expansion of $p_2({\ell}^*(z))$ around   $z=0$ is given by
$\frac{\bar{p}_{2,1}(z)}{z^{k(m-n)}}\,G(z)$, where $G(z)$ is the
Taylor series of $1/\bar{p}(z)$ at $z=0$. Observe that $G(z)$
exists since all the derivatives of $1/\bar{p}(z)$ at $z=0$ exist
(note that the denominator of all the derivatives is a power of the
polynomial $\bar{p}(z)$, and $\bar{p}(0)=d_n\nu(0)^n=d_n
b_k^n\not=0$).  In addition, taking into account that
$$\nu^{j)}(0)=b_{k-j},\,\,0\leq j\leq k,\quad \mbox{and}\quad
\nu^{j)}(0)=a_{j-k},\,j\geq k+1,$$ and that $\frac{\partial^j
(1/\bar{p}(z))}{\partial z^j}_{|z=0}$ is obtained from
$\nu^{i)}(0),\,0\leq i\leq j$, we get that
\[G(z)=\frac{1}{\bar{p}(0)}+z\frac{\partial (1/\bar{p}(z))}{\partial
z}_{|z=0}+\cdots= h_0(b_k)+\cdots+z^k
h_k(b_k,\ldots,b_{0})\]\[+z^{k+1}
h_{k+1}(b_k,\ldots,b_{0},a_{1})+\cdots+z^{k+u}
h_{k+u}(b_k,\ldots,b_{0},a_{1},\ldots,a_{u})+ \cdots,\]
where $h_j(b_k,\ldots,b_{0},a_{1},\ldots,a_{j-k}),\,j\geq 0,$ denotes a rational function depending on $b_k,\ldots,b_{0},a_{1},\ldots,a_{j-k}$.

\para

As we stated above, we need to determine the terms with non positive
exponent in
$$p_2({\ell}^*(z))=\frac{\bar{p}_{2,1}(z)}{z^{k(m-n)}}\,G(z).$$
In the following, we prove that they can be obtained by just
computing $b_k,\ldots,b_{0},a_{1},\ldots, a_{km}$. Indeed:
\begin{enumerate}
\item[1.1.] Let $m=n$. Then,  we need to compute  the terms with non positive exponent  in
$$\bar{p}_{2,1}(z)G(z)=(c_m\nu^{m}+c_{m-1}z^{k}\nu^{m-1}+c_{m-2}z^{2k}\nu^{m-2}+\cdots+c_{0}z^{km})$$
$$(h_0(b_k)+\cdots+z^k h_k(b_k,\ldots,b_{0})+z^{k+1} h_{k+1}(b_k,\ldots,b_{0},a_{1})+\cdots).$$
Thus, we only need the independent term   $c_m
b_{k}^{m}h_0(b_k)$.
\item[1.2.] Let $m<n$. In this case, we need to determine  the terms with non positive exponent  in $z^{k(n-m)}\bar{p}_{2,1}(z)G(z).$
However, since $n-m>0$, we conclude that there are no such terms.
   \item[1.3.] Let $m>n$. Then, we need to compute  the terms with non positive exponent
   in ${\bar{p}_{2,1}G}/{z^{k(m-n)}}$ which implies that we need to determine the terms having degree less or equal to  $k(m-n)$
   in the product $\bar{p}_{2,1}(z)G(z)$. Those terms are included in the product\\

\noindent
$(c_m(b_0 z^{{k}}+b_{1}z^{k-1}+\cdots
+b_{k-1}z+b_{k})^{m}+c_{m-1}z^{k}(b_0
z^{{k}}+b_{1}z^{k-1}+\cdots +b_{k-1}z+b_{k})^{m-1}+\cdots+c_0z^{k
m})\cdot
(h_0(b_k)+\cdots+z^k h_k(b_k,\ldots,b_{0})+z^{k+1} h_{k+1}(b_k,\ldots,b_{0},a_1)+\cdots +z^{k(m-n)}h_{k(m-n)}(b_k,\ldots,b_{0},a_{1},\ldots,a_{k(m-n)}))$\\

(we  do not include the term
$z^kB(z^N)$ in this product since after multiplying,  it only provides terms of degree
greater than $km$). Therefore, at most we have to compute $\ell(z)$
till the terms $b_k, \ldots, b_{0}, a_{1},\ldots,a_{k(m-n)}$ appear. That is, $k+1+k(m-n)$ terms are needed.
    \end{enumerate}

Taking into account the cases $1.1, 1.2$, and $1.3$, we deduce that
at most we have to compute $k+1+k(m-n)$ terms in $\ell(z)$. Finally,
we prove that $k+1+k(m-n)\leq 2\deg(p_1)+1$. For this purpose, let
$d(r_2)$ denote the maximum exponent of $z$ in $r_2(z)$. We observe
that $d(r_2)\leq 1$; otherwise, since
$$F(z:r_2(z):r_3(z):1)=F(z/r_2(z):1:r_3(z)/r_2(z):1/r_2(z))=0 $$
(for $|z|>M$) by continuity, we get
$$\lim_{z\rightarrow\infty}F(z/r_2(z):1:r_3(z)/r_2(z):1/r_2(z))=F(0:1:C:0)=0$$
where $C:=\lim_{z\rightarrow\infty}r_3(z)/r_2(z)$. If $C\in {\Bbb C}$, we get that $(0:1:C:0)$ is an infinity point of the input curve which is impossible
since we have assumed that the input curve does not have infinity points of the form $(0:a:b:0)$. If $C=\infty$, we reason as above but we divide by
$r_3(z)$. In this case, we get the infinity point $(0:0:1:0)$ which is again impossible.

On the other hand, since
$r_2(z)=p_2(\ell(z^{-1}))=\frac{p_{21}(\ell(z^{-1}))}{p(\ell(z^{-1}))}$,
we get that $d(r_2)=(m-n)k/N$, where $m=\deg(p_{21})$ and
$n=\deg(p)$ (see Chapter 4 in \cite{Walker}). Hence, $(m-n)k/N\leq 1$
which implies that $(m-n)k\leq N$. In addition, since
$N\leq\deg_s(p(s)-tp_{11}(s))=\deg(p_1)$ (see Remark 4 in
\cite{BlascoPerezII}), we get that $k+1+k(m-n)\leq 2k(m-n)+1\leq
2\deg(p_1)+1$.

 \item  Let us assume that $b_k=0$ for $k>0$. That is, there are no terms with negative exponent in $\ell(z)$. Then, we write $\ell(z):=b_0+B(z)$, where
 \[B(z)=\sum_{j=1}^\infty
a_{j}z^{q_{j}/N},\,\,N
\in\mathbb{N}^+,\,q_{j}\in\mathbb{N}^+,\,0<q_1<q_2<\cdots,\,
a_{j}\in\mathbb{C}\setminus{\{0\}},\] and
 \[{\ell}^*(z):=\ell(z^N)=b_0 +B(z^N)=b_0+z^{q_1}\left(a_1+\sum_{j=2}^\infty a_{j}z^{q_{j}-q_1}\right),\,\, B(z^N)=\sum_{j=1}^\infty
a_{j}z^{q_{j}}.\] In this case, we denote $\nu:=b_0+z^{q_1}(a_1+\sum_{j=2}^\infty
a_{j}z^{q_{j}-q_1})$. In addition, we write
$$p(t)=p^*(t)(t-b_0)^r,\quad \gcd(p^*(t),t-b_0)=1\quad \mbox{for
some}\, \, r\in {\Bbb N}.$$ Under these conditions, we get that
$p_2({\ell}^*(z))=$
 \[\frac{p_{2,1}(\nu)}{p(\nu)}=\frac{p_{2,1}(\nu)}{p^*(\nu)(\nu-b_0)^r}=\frac{p_{2,1}(\nu)}{z^{rq_1}p^*(\nu)(a_1+\sum_{j=2}^\infty a_{j}z^{q_{j}-q_1})^r}:= \frac{\bar{p}_{2,1}(z)}{z^{rq_1}\bar{p}(z)},\]
 where
 \[\bar{p}_{2,1}(z)=p_{2,1}(\nu)=c_m\nu^{m}+c_{m-1}\nu^{m-1}+\cdots+c_{0},\,c_m\not=0,\,m=\deg(p_{2,1})\]
 and $\bar{p}(z)=p^*(\nu)(a_1+\sum_{j=2}^\infty a_{j}z^{q_{j}-q_1})^r=$
\[(d_n\nu^{n}+d_{n-1}\nu^{n-1}+\cdots+d_{0})\left(a_1+\sum_{j=2}^\infty a_{j}z^{q_{j}-q_1}\right)^r,\,d_n\not=0,\,n:=\deg(p^*).\]
The generalized series expansion of $p_2({\ell}^*(z))$ around $z=0$
is given by $\frac{\bar{p}_{2,1}(z)}{z^{rq_1}}\,G(z)$, where $G(z)$
is the Taylor series of $1/\bar{p}(z)$ at $z=0$. Observe that $G(z)$
exists since all the derivatives of $1/\bar{p}(z)$ at $z=0$ exist
(note that the denominator of all the derivatives is a power of the
polynomial $\bar{p}(z)$, and $\bar{p}(0)=p^*(\nu(0))a_1=
p^*(b_0)a_1\not=0$). Reasoning as in case $1$, one may check that
$G(z)=\frac{1}{\bar{p}(0)}+z\frac{\partial (1/\bar{p}(z))}{\partial
z}_{|z=0}+\cdots=$ \[= h_0(b_0, a_1)+z h_1(b_0,a_1,
a_{2})+\cdots+z^k h_k(b_0,a_1,\ldots,a_{k+1})+ \cdots,\] where
$h_j(b_0,a_1,\ldots,a_{j+1}),\,j\geq 0$ is a rational function
depending on $b_0,a_1,\ldots,a_{j+1}$. \para

Since we need to compute  the terms with non positive exponent  in
$$p_2({\ell}^*(z))=\frac{\bar{p}_{2,1}(z)}{z^{rq_1}}\,G(z),$$ we reason
as in case $1.1$ (if $r=0$), or case $1.3$ (if $r>0$), and we
conclude that at most, we have to determine $\ell(z)$ till the terms
$b_0, a_1,\ldots,a_{rq_1+1}$ appear. That is,  in this case, at most
$rq_1+2$ terms are needed. Finally, we prove that $rq_1+2\leq
2\deg(p_1)+1$. For this purpose,  we reason as above and since
$$r_2(z)=\frac{p_{21}(\ell(z^{-1}))}{p(\ell(z^{-1}))}=\frac{p_{21}(\ell(z^{-1}))}{(\sum_{j=1}^\infty a_{j}z^{-q_{j}/N})^rp^*(\ell(z^{-1}))},$$
and
$\lim_{z\rightarrow\infty}p_{21}(\ell(z^{-1}))/p^*(\ell(z^{-1}))=p_{21}(b_0)/p^*(b_0)\in
{\Bbb C}$ (and thus, $d(p_{21}(\ell(z^{-1})))=d(p^*(\ell(z^{-1})))$, we get that  $d(r_2)=rq_1/N$ (see Chapter 4 in \cite{Walker}). Since $d(r_2)\leq 1$, we deduce that $rq_1\leq N$. In
addition, since $N\leq\deg(p_1)$ (see Remark 4 in
\cite{BlascoPerezII}), we get that $rq_1\leq \deg(p_1)$, and thus
$rq_1+2\leq  \deg(p_1)+2\leq 2\deg(p_1)+1$. \hfill
$\Box$\end{enumerate}

\end{document}